\documentclass[twoside,11pt]{amsart}
\usepackage{graphicx, lscape}
\usepackage{amsmath, amsfonts, amssymb, ifthen}

\renewcommand{\Re}{\operatorname{Re}}
\renewcommand{\Im}{\operatorname{Im}}
\renewcommand{\(}{\left(}
\renewcommand{\)}{\right)}  
\newcommand{\abs}[1]{\left\arrowvert #1\right\arrowvert}

\newcommand{\field}[1]{\mathbb{#1}}
\newcommand{\CC}{\ensuremath{\field{C}}} 
\newcommand{\RR}{\ensuremath{\field{R}}} 
\newcommand{\NN}{\ensuremath{\field{N}}} 
 
\newcommand{\ZZ}{\ensuremath{\field{Z}}} 
\newcommand{\Ct}{\ensuremath{\field{C}^2}}
\newcommand{\Ch}{\ensuremath{\hat{\field{C}}}}
 
\newcommand{\Pt}{\ensuremath{\field{P}^2}}

\newcommand{\Nvs}{\mathcal{N}_{v}^{\#}}
 
\newcommand{\eps}{\epsilon}
\newcommand{\lam}{\lambda}

\newcommand{\A}{A(C_{J_p})}
\newcommand{\Apt}{A_{pt}(C_{J_p})}
\newcommand{\Acc}{A_{cc} (C_{J_p})} 

\theoremstyle{plain}
\newtheorem{thm}{Theorem}[section]
\newtheorem{cor}[thm]{Corollary}
\newtheorem{lem}[thm]{Lemma}
\newtheorem{prop}[thm]{Proposition}

\newtheorem{example}[thm]{Example}

\theoremstyle{definition}

\newtheorem{rem}[thm]{Remark}
\newtheorem{question}[thm]{Question}

\newcommand{\drawfigtwbascl}{\scalebox{.5}{\includegraphics{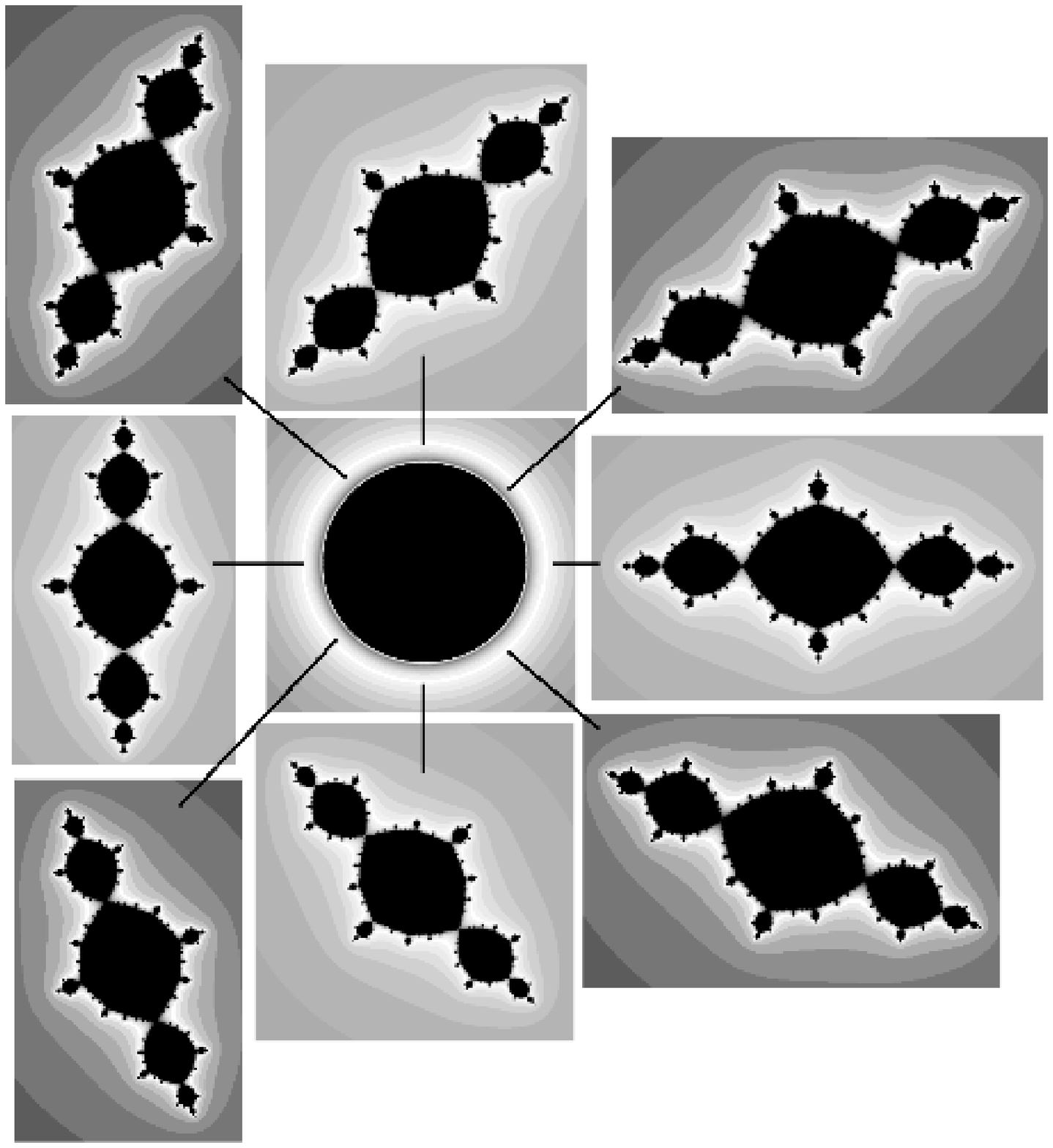}}}

\newcommand{\drawfigaerocant}{\scalebox{.52}{\includegraphics{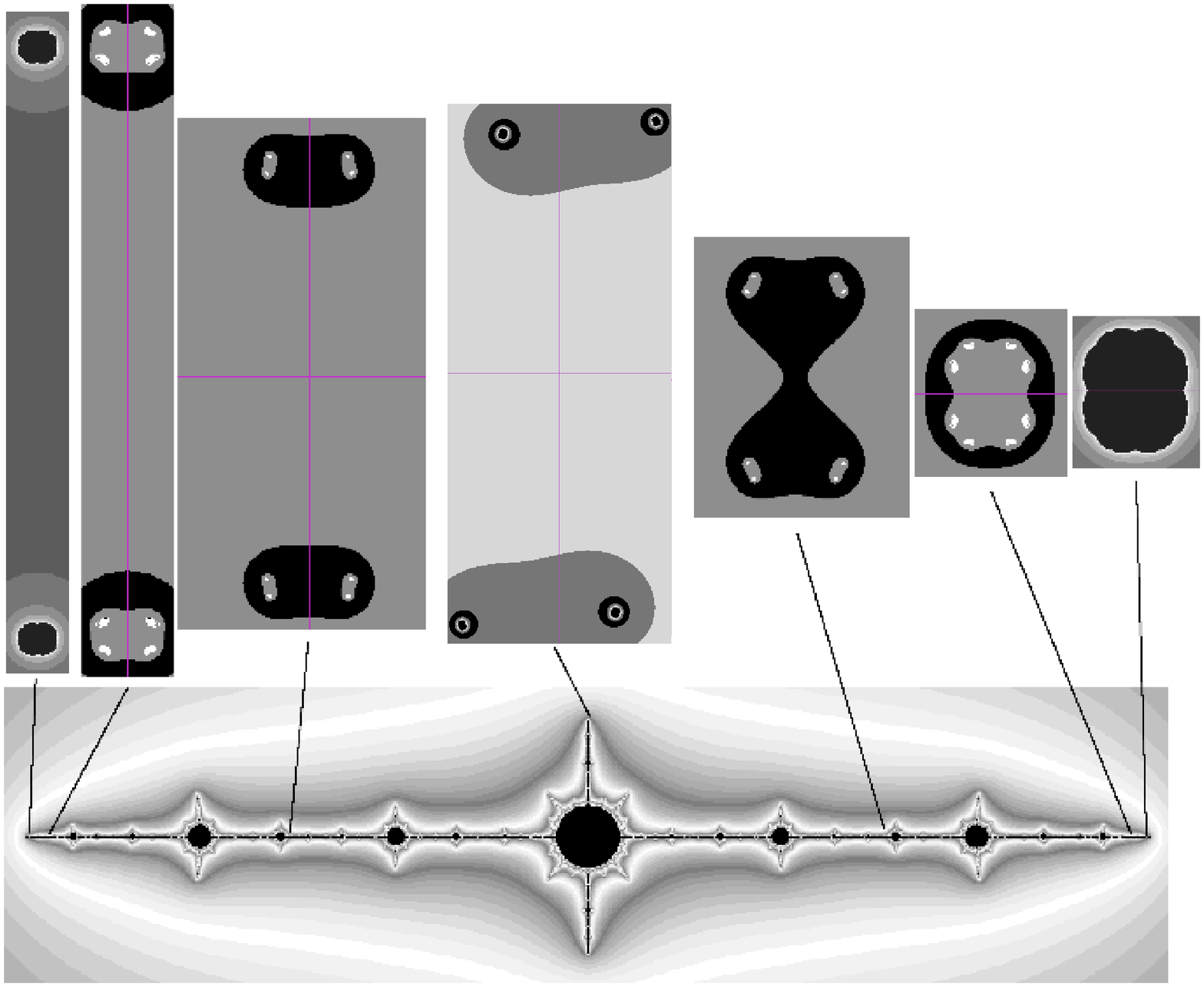}}}

\newcommand{\drawfigcantcircbas}{\scalebox{.8}{\includegraphics{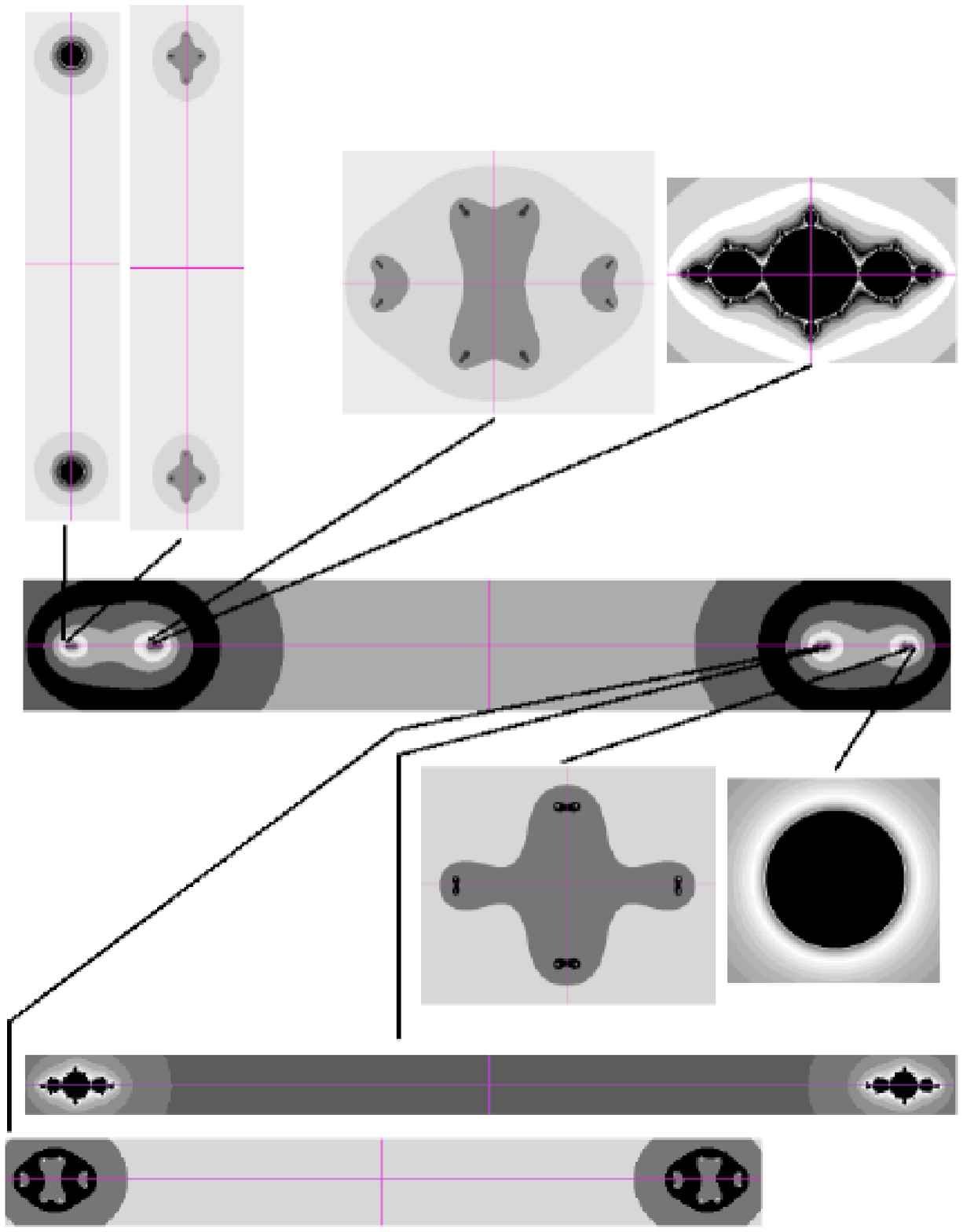}}}

\newcommand{\showcomments}{yes}

\newsavebox{\commentbox}
\newenvironment{comment}%
{\ifthenelse{\equal{\showcomments}{yes}}%
{\footnotemark
    \begin{lrbox}{\commentbox}
    \begin{minipage}[t]{1.25in}\raggedright\sffamily\small
    \footnotemark[\arabic{footnote}]}
{\begin{lrbox}{\commentbox}}}%
{\ifthenelse{\equal{\showcomments}{yes}}%
{\end{minipage}\end{lrbox}\marginpar{\usebox{\commentbox}}}
{\end{lrbox}}}

\begin{document}

\title{Axiom A Polynomial skew products of $\Ct$ \\ and their postcritical sets}

\author[L. ~DeMarco]{Laura DeMarco$^1$}
\address{Department of Mathematics\\
University of Chicago\\
5734 S. University Ave.\\
Chicago, IL 60637\\
USA}
\email{demarco@math.uchicago.edu}
\author[S.L. ~Hruska]{Suzanne Lynch Hruska$^2$}
\address{Department of Mathematical Sciences\\
University of Wisconsin Milwaukee\\
PO Box 413\\
Milwaukee, WI 53201\\
USA}
\email{shruska@msm.umr.edu}

\date{\today}

\begin{abstract}
A polynomial skew product of $\Ct$ is a map of the form $f(z,w) = (p(z), q(z,w))$, where $p$ and $q$ are polynomials, such that $f$ extends holomorphically to an endomorphism of $\mathbb{P}^2$ of degree $\geq 2$.
For polynomial maps of $\CC$, hyperbolicity is equivalent to the condition that the closure of the postcritical set is disjoint from the Julia set; further,
critical points either iterate to an attracting cycle or infinity. 
For polynomial skew products, Jonsson (\cite{MatSk}) established that $f$ is Axiom A if and only if 
 the closure of the postcritical set is disjoint from the right analog of the Julia set.  Here we present an analogous conclusion:  critical orbits either escape to infinity or accumulate on an attracting set.    
In addition, we construct new examples of Axiom A maps demonstrating various postcritical behaviors.
\end{abstract}

\maketitle

\markboth{\textsc{L. DeMarco and S. Hruska}}
  {\textit{Axiom A Polynomial Skew Products}}

\footnotetext[1]{$^{,2}$Research supported in part by a grant from
the National Science Foundation.}


\section{Introduction}
\label{sec:introduction}

A polynomial skew product of $\Ct$ is a map $f(z,w) = (p(z), q(z,w))$, where $p$ and $q$ are polynomials.  We will assume that $f$ is regular of degree $d\geq 2$, meaning that it extends holomorphically to an endomorphism of $\mathbb{P}^2$;  consequently the polynomials $p$ and $q$ both have degree $d$.  

In this paper, we study the postcritical set for Axiom A polynomial skew products, and we construct new examples of Axiom A maps demonstrating various postcritical behaviors.  This is motivated by results in one-dimensional complex dynamics relating the behavior of the postcritical set to hyperbolicity.

For polynomial maps of $\CC$, hyperbolicity is equivalent to the condition that the postcritical set is disjoint from the Julia set; by the classification of Fatou components, critical points are either in the basin of an attracting cycle or they escape to infinity under iteration (see e.g.\ \cite[Theorem 3.13]{McMullen:CDR}).  
For polynomial skew products, Jonsson established that Axiom A is characterized by postcritical behavior
\cite[Cor 8.3]{MatSk}:  the postcritical set must be disjoint from the right analog of the Julia set.  Here we present the analogous conclusion:  critical points either escape to infinity or are in the basin of an attracting set.

\medskip\noindent{\bf Dynamics over $J_p$.}  
Let $J_p$ denote the Julia set of the base polynomial $p$.  We are primarily
interested in the dynamics of $f$ restricted to $J_p\times\CC$.  This invariant
subset contains $J_2$, the closure of the repelling cycles for $f$, which coincides
with the support of the measure of maximal entropy for $f$ \cite{MatSk}.  

When $f$ is Axiom A, its nonwandering set is contained in $(J_p \cup A_p) \times \CC$, 
where $A_p$ is the finite set of attracting periodic points of $p$.  
The dynamics of $f$ over $A_p$ reduces to one-dimensional complex dynamics, 
so the nontrivial part of the saddle set of $f$ is also contained in $J_p\times\CC$.  

As $f$ preserves the family of vertical lines $\{ z \} \times \CC$, 
the iterates $f^n|_{ \{z \} \times \CC}$ form a composition of polynomial maps of 
$\CC$ of degree $d$, with bounded coefficients for $z\in J_p$.  
The dynamics of such composition sequences have been studied
in \cite{FSrand} and \cite{MCrand}.  In our 
study of skew products in $\Ct$, we combine results about (hyperbolic) 
composition sequences with one- and two-dimensional complex iteration theory.  We remark
that Sumi also studies a notion of hyperbolicity for composition sequences in the setting of rational skew products  \cite{sumi:etds01, sumi:etds06}. 

\medskip\noindent{\bf Main result.} 
If $f(z,w) = (p(z), q(z,w))$ is Axiom A, 
let $\Lambda$ denote the union of the basic sets of saddle type in $J_p \times \CC$
(see \S\ref{sec:background} for definitions). 
Its unstable manifold $W^u(\Lambda)$ consists of all points $x\in\CC^2$ for which there exists a backward orbit $x_{-k}$ of $x$ converging to $\Lambda$.  
Axiom A implies  $W^u(\Lambda) \cap J_2 = \emptyset$.
Our main theorem shows that in $J_p\times\CC$, each of $\Lambda$ and $W^u(\Lambda)$ plays the role of the one-dimensional attracting periodic orbits, in the following sense.

Set $q_z(w) = q(z,w)$, and define the {\bf critical locus over $J_p$} by 
 $$C_{J_p} = \{ (z,c) \colon z \in J_p, q'_z(c)=0 \}.$$  
If $X$ is any subset of $\Ct$, its {\bf accumulation set}\begin{footnote}{For a point $x \in \Ct$,  $A(x)$ coincides with the $\omega$-limit set of $x$.}\end{footnote}  is 
 $$A(X) = \displaystyle{\bigcap_{N \geq 0} \  \overline{ \bigcup_{n\geq N}  f^n(X) }}.$$
We will study the accumulation set $A(C_{J_p})$ as well as the 
pointwise and component-wise accumulation sets,
 $$\Apt = \overline{\bigcup_{x \in C_{J_p}} A(x)} \quad\mbox{ and } \quad
 \Acc = \overline{\bigcup_{C \in \mathcal{C}(C_{J_p})} A(C) }, $$
where $\mathcal{C}(C_{J_p})$ denotes the collection of connected components of $C_{J_p}$.
Each of these sets is closed, and we clearly have $\Apt \subseteq \Acc \subseteq \A$.

\begin{thm} \label{thm:Achain}
If $f$ is an Axiom A polynomial skew product of $\Ct$, then 
$$\Lambda = \Apt \subseteq \Acc \subseteq \A = W^u(\Lambda) \cap (J_p \times \CC).$$  
\end{thm}

In particular, the first equality $\Lambda = \Apt$ shows that each critical point in $J_p\times\CC$
either tends to the saddle set or escapes to infinity.  
To explain the appearance of $W^u(\Lambda)$, we can compare to the invertible setting.  The real, or smooth, $\lambda$-lemma states that  if $p$ is a saddle periodic point, then the 
forward images of a disk transverse to $W^s(p)$ which intersects a neighborhood of $p$ 
will tend to $W^u(p)$  (see e.g.\ \cite[Chp.\ 5, Theorem 11.1]{Rob}).  
In our setting, the critical locus $C_{J_p}$ (in fact any set transverse to $W^s(\Lambda)$ 
and disjoint from $J_2$) is like such a disk (see Proposition \ref{prop:Wattracts}).

\medskip\noindent{\bf Explicit families of Axiom A examples.}
Among dynamical systems with chaotic behavior, Axiom A maps are the most tractable.  
Jonsson has shown that the Axiom A maps form an open subset of the parameter space 
of all skew products, which allows him to define \textbf{hyperbolic components} as connected components of the subset of Axiom A maps (\cite[Corollary 8.15, Definition 8.16]{MatSk}).

In this paper, we construct examples of Axiom A maps supporting different chains of  equalities or inequalities in Theorem~\ref{thm:Achain}, and distinguish the hyperbolic components of our examples.    This work yields the next two results.

\begin{prop}  \label{prop:products}
Polynomial skew products in the same hyperbolic component as a product have equality for every 
inclusion of Theorem \ref{thm:Achain}.
\end{prop}

\begin{thm} \label{thm:examples}
There exist examples of Axiom A polynomial skew products 
$f_1, f_2, f_3, f_4$, none of which are in the hyperbolic component of a product, such that:
\begin{enumerate}
\item for $f_1$, $\Apt = \Acc = \A = \emptyset$;
\item for $f_2$, $\Apt = \Acc = \A \not= \emptyset$;
\item for $f_3$, $\Apt \neq \Acc = \A$;
\item for $f_4$, $\Apt = \Acc \neq \A$.
\end{enumerate}
Moreover, each $f_k$ is also Axiom A as an endomorphism from $\mathbb{P}^2$ to itself.
\end{thm}

There are few known examples of Axiom A endomorphisms in two (or more) dimensions.   
Fornaess and Sibony study examples of Axiom A endomorphisms of $\Pt$ in \cite{FSexamples, FShyp}.
Among polynomial skew products, Jonsson (\cite[\S 9]{MatSk}) lists known Axiom A maps to be: products of hyperbolic polynomial maps of $\CC$, small perturbations of hyperbolic products, and he gives a 
degree 2 example of an Axiom A map not in the same hyperbolic component as any product.
Diller and Jonsson (\cite{MatDil}) generalize this last example for all $d \geq 2$. Their examples also satisfy (4) in Theorem \ref{thm:examples}.   Sumi has communicated to the authors another new example of an Axiom A polynomial skew product satisfying (2) of Theorem~\ref{thm:examples}, which is very different from ours.  We describe his construction in Example~\ref{exmp:sumi}.

\medskip\noindent{\bf Stability and holomorphic motions.} 
One of our main tools for distinguishing hyperbolic components 
is holomorphic motions of $J_2$.
Jonsson established a general result about holomorphic motions of hyperbolic sets of holomorphic endomorphisms  \cite[Theorem C]{MatHM}, which implies that if $f$ is an Axiom A polynomial skew product, then $J_2$ moves holomorphically under perturbation. Further, we show: 

\begin{thm}  \label{thm:motion}
For a holomorphic family of Axiom A polynomial skew products, the holomorphic motion
of $J_2$ preserves the vertical fibration, inducing a holomorphic motion of each 
fiber Julia set $J_z = J_2 \cap (\{z\}\times\CC)$. 
\end{thm}

A more precise statement is given in Theorem \ref{thm:SkewHM}.
As a corollary, if a polynomial skew product is in the same hyperbolic component as a product, then all fiber Julia sets $J_z$ are homeomorphic 
(Corollary \ref{cor:productcomponent}). 
This answers a question from \cite{SLHSK}, originally posed to the author by Eric Bedford.
Our examples satisfying (1) and (2) of Theorem~\ref{thm:examples} will show that the converse is not true, for in those examples all fiber Julia sets are homeomorphic, but the maps are not in the same hyperbolic component as any product.

Theorem \ref{thm:motion} also implies that the equalities $\Apt = \Acc$ and $\Apt = \A$
are preserved under perturbation (Propositions \ref{prop:Lcomponent} and \ref{prop:AptisAcc}).  
We leave it as an 
open question to show that the equality $\Acc=\A$ is also preserved (see Question \ref{question2}).

\medskip\noindent{\bf Organization of sections.} 
In Section~\ref{sec:background} we provide some needed prerequisite material on the dynamics of polynomial skew products, hyperbolic maps of $\Pt$, and general composition sequences of polynomials.

In Section~\ref{sec:ACS1}, we prove Theorem~\ref{thm:Achain}.  

In Section~\ref{sec:fiberhomeo}, we prove Theorem~\ref{thm:motion} 
on holomorphic motions.  

In Section~\ref{sec:L=W}, we establish conditions giving equality for all inclusions in Theorem~\ref{thm:Achain}, prove Proposition~\ref{prop:products}, and  provide specific examples of maps satisfying (1) and (2) of Theorem~\ref{thm:examples}.

In Section~\ref{sec:aerocantor}, we establish a necessary and sufficient condition for $\Apt = \Acc$ in hyperbolic components (Lemma~\ref{lem:AptAccchar}), and construct Axiom A skew products  with base Julia set connected, but topologically varying fiber Julia sets (some connected, some disconnected).  We show these maps satisfy (3) of Theorem~\ref{thm:examples}.   These maps are a perturbation of \cite[Example 9.7]{MatSk} (which is not Axiom A).

In Section~\ref{sec:matnonprodgen}, we produce Axiom A maps illustrating (4) of Theorem \ref{thm:examples}.  These are a generalization of \cite[Example 3.9]{MatDil}.  

We close the paper with some open questions.

\subsection*{Acknowledgements}  We  thank Eric Bedford and Curt McMullen 
for helpful conversations during our work on this project, 
Greg Buzzard for feedback on drafts of the paper, Mattias Jonsson for pointing
out Remark~\ref{rem:semiconj}, Lois Kailhofer for providing insight on natural extensions and inverse limits, and Hiroki Sumi for describing Example~\ref{exmp:sumi}.  We would also like to thank the referee for his careful reading of the paper and many constructive suggestions.  

\section{Background}
\label{sec:background}

A polynomial map $f: \Ct\to\Ct$ is a (regular) {\bf polynomial skew product} if it has the form
  $$f(z,w) = (p(z), q(z,w)),$$
for polynomials $p$ and $q$ of degree $>1$, and if $f$ extends analytically as a map
from $\Pt$ to itself.  We refer the reader to \cite{MatSk} for a general treatment of 
polynomial skew products and proofs of many of the statements below.  

\medskip\noindent\textbf{Global dynamics}.
Fix a polynomial skew product $f(z,w) = (p(z), q(z,w))$.  
Let $K_p$ denote the filled Julia set of $p$ (the set of points with bounded orbit) 
and $J_p = \partial K_p$ its Julia set.  
Let $K = \{x\in\Ct: \sup_n \|f^n(x)\| < \infty\}$ be the filled Julia set 
of $f$.  Set $K_z = K\cap (\{z\}\times\CC)$ and $J_z = \partial K_z$ in $\{z\}\times\CC$.   The set 
	$$J_2 = \overline{ \bigcup_{z\in J_p} J_z }$$
is the closure of the repelling cycles of $f$ and also the support of its measure of 
maximal entropy.   The {\bf nonwandering set} $\Omega$ of $f$ is the set of points having no neighborhood $V$
such that $f^n(V)\cap V = \emptyset$ for all $n>0$.

\medskip\noindent\textbf{Vertical dynamics and expansion}.
For a fixed $z\in K_p$, we define $q_z(w) = q(z,w)$, and set 
	$$Q^n_z = q_{p^{n-1}(z)}\circ \cdots\circ q_{p(z)}\circ q_z.$$
We have $q_z K_z = K_{p(z)}$ and $q_z J_z = J_{p(z)}$.   
Let $C_z$ be the critical set of $q_z (w)$ in $\{z\}\times\CC$.

Let $Z \subseteq K_p$ be compact with $p(Z) \subseteq Z$ (e.g.\ 
$Z=J_p$ or $Z$ is an attracting periodic point of $p$).   Set 
 $$C_Z = \bigcup_{z \in Z} C_z, \qquad PC_Z = \overline{ \bigcup_{n\geq 1} f^n C_Z  },$$
and
 $$J_Z = \overline{ \bigcup_{z \in Z} J_z}.$$

We say $f$ is {\bf vertically expanding over $Z$} if there exist constants 
$c>0$ and $L >1$ such that 
  $$\abs{(Q^n_z)'(w)} \geq cL^n$$ 
for all $z\in Z, w\in J_z$, and all $n\geq 1$.

In this paper, we repeatedly use the following results.

\begin{prop}  \cite[Proposition 2.1]{MatSk}  \label{prop:Jprop2.1}
For every polynomial skew product $f$, 
$z \mapsto K_z$ is upper semi-continuous (in the Hausdorff topology), 
$z \mapsto J_z$ is lower semi-continuous, and 
the sequence of polynomials
$\{ Q_z^n \}_{n\geq 1}$ is normal exactly on $\{ z \} \times (\CC \setminus J_z)$.
\end{prop}

\begin{prop} \cite[Proposition 2.3]{MatSk} \label{prop:Jprop2.3}
For any $z \in \CC$, $J_z$ is connected if and only if 
$C_{p^n(z)} \subset K_{p^n(z)}$ for all $n \geq 0$.
\end{prop}

\begin{thm}  \cite[Proposition 3.5]{MatSk}  \label{thm:Jprop3.5}
If $f$ is vertically expanding over $Z$, then $z\mapsto J_z$ is continuous 
in the Hausdorff topology with $z\in Z$.  In particular, if $f$ is vertically 
expanding over $J_p$, we have $$J_2 = \bigcup_{z\in J_p} J_z.$$
\end{thm}

\begin{thm}  \cite[Theorem 3.1]{MatSk}  \label{thm:Jthm3.1}
$f$ is vertically expanding over $Z$ if and only if $PC_Z \cap J_Z = \emptyset$.
\end{thm}

We do not define Axiom A in this paper, because we use only the following:

\begin{thm}  \cite[Theorem 8.2]{MatSk}  \label{thm:Jthm8.2}
A polynomial skew product $f$ is Axiom A on $\Ct$ if and only if 
\begin{itemize}
\item[(i)] $p$ is expanding on $J_p$ (i.e.~$p$ is hyperbolic);
\item[(ii)] $f$ is vertically expanding over $J_p$; and
\item[(iii)] $f$ is vertically expanding over $A_p$.  
\end{itemize}
Moreover, if $f$ is Axiom A on $\Ct$, then the nonwandering set equals the closure of the set of periodic points of $f$ (equals the chain recurrent set).  

Further, $f$ is Axiom A on $\Pt$ if in addition, 
\begin{itemize}
\item[(iv)]  the extension of $f$ to the line at infinity of $\Pt$ is hyperbolic as a one-dimensional polynomial.
\end{itemize}
\end{thm}

\medskip\noindent\textbf{Hyperbolicity for endomorphisms of $\Pt$}.
For a non-invertible holomorphic mapping $f$ of $\Pt$, basic results and definitions for hyperbolicity and stability of an invariant set $L$  must be given in terms of the {\bf natural extension}, $\hat{L}$, which is the space of all sequences of prehistories (backward orbits) in $L$, with $f|_L$ inducing $\hat{f}|_{\hat{L}}$ which is a shift.  See, for example, \cite[\S A.2]{MatSk}, \cite[\S2]{FShyp}, and \cite[\S1.1]{MatDil}
for definitions and properties.

 If $f$ is Axiom A, then its nonwandering 
set $\Omega$ decomposes into a finite union of basic sets $\bigcup \Omega_i$ such that $f(\Omega_i)
=\Omega_i$ and $f$ is transitive on each $\Omega_i$.  A basic set $\Omega_i$ 
is of saddle type if its unstable (complex) dimension is 1.  

For an Axiom A polynomial skew product, we let ${\Lambda}$ denote the union of basic saddle sets 
contained in $J_p\times\CC$.  On these sets, $f$ is contracting in the fiber direction.   
Theorem~\ref{thm:Jthm8.2} also yields that $\Lambda$ is the closure of the set of saddle periodic points $(z,w)$ with $z \in J_p$, since this theorem establishes that periodic points are dense in the nonwandering set.

The set $J_2$ is the subset of $\Omega$ of unstable
dimension 2.  

The {\bf stable manifold} $W^s(L)$ of an invariant set $L$ is the set of all points $x\in \Ct$ for which 
the forward orbit $\{ f^k(x) \}_{k \geq 0}$ converges to $L$.    
The {\bf unstable manifold} $W^u(L)$ is the set of all points $x\in \Ct$ for which 
there exists a prehistory (backward orbit) $\{ x_{-k} \}_{k \geq 0}$, converging to $L$, i.e., $x = x_0, f(x_{-k}) = x_{-k+1}$ for all $k \geq 0$,  and $x_{-k} \to L$ as $k \to \infty$.   

In \cite[Corollary 8.14]{MatSk}, Jonsson states that the natural extension of the nonwandering set of an Axiom  A polynomial skew product is stable;  i.e., if $g$ is a holomorphic map of $\Pt$ which is $C^1$-close to $f$, then there is a homeomorphism $\Phi \colon \hat{\Omega}_f \to \hat{\Omega}_g$ conjugating $\hat{f}$ to $\hat{g}$ (and $\Phi$ respects decomposition into basic sets), and $\Phi$ can be chosen close to the identity.  In fact, by \cite[Theorem B]{MatHM}, $\hat{\Omega}$ moves holomorphically, in a sense which we will not define precisely.
\begin{footnote}{In this paper, the only stability result we use is that $J_2$ moves holomorphically for Axiom A polynomial skew products (see Section \ref{sec:fiberhomeo}), but we state this alternate result for benefit of the reader.} \end{footnote}

\medskip\noindent\textbf{General composition sequences of polynomial maps of $\CC$}.
Fornaess and Sibony (\cite{FSrand}) study fundamental
properties of Julia sets arising from compositions of sequences of arbitrarily chosen polynomial maps of $\CC$, with uniformly bounded degrees and coefficients.  
Comerford (\cite{MCrand}) establishes an analog of the one-variable postcritical characterizations of hyperbolicity
in this setting. We  apply Comerford's results to the composition of fiber maps $q_z$ along the orbit of a point $z \in J_p$.   We state a key result here for our setting of skew products, though it was written for
general composition sequences (with uniform bounds on  degree and coefficients).

\begin{thm}  \cite[Theorem 4.1]{MCrand} \label{thm:Cthm4.1}  
Suppose $f$ is a polynomial skew product which is vertically expanding over $J_p$.
Let $T$ be a compact, connected subset of $\{ z \} \times (\CC \setminus J_z)$, for some $z\in J_p$.  
There exist $\lambda \in (0,1)$ and $\kappa > 0$ such that for all $m \in \NN$, 
$$\text{diam}^{\#}(f^m(T)) \leq \kappa \lambda^m \text{ diam}^{\#}(T),$$
where diam$^{\#}$ is the diameter in the spherical metric on $\Ch$.  
Further, $\lambda$ depends only on  $f$, and $\kappa$ in addition 
depends on the distance from $T$ to $J_z$ (in the spherical metric).  
\end{thm}

\section{Saddle sets and the attractor of the postcritical set}
\label{sec:ACS1}

Throughout this section, let $f$ be an Axiom A polynomial skew product of degree $d$.  Recall that $\Lambda$ denotes the union of the saddle basic sets in $J_p \times \CC$, and is the closure of the saddle periodic points $(z,w)$ with $z \in J_p$.  In this section we prove Theorem~\ref{thm:Achain}. 


\bigskip\noindent{\bf Trapping radius.}
If $w \in \CC$, let $B^{\#}(w,r)$ denote the disk about $w$ of radius $r$ in the spherical metric on $\Ch$.
For $S \subset J_p \times \CC$, the {\bf vertical neighborhood} of $S$ of radius $r$ is the set 
 $$\Nvs(S, r) = \bigcup_{(z,w)\in S} \{z\} \times B^{\#}(w,r).  $$
The following lemma is a corollary of Theorem \ref{thm:Cthm4.1}, and was
observed by Comerford in the setting of general composition sequences.

\begin{lem} \label{lem:trap}
Let $T$ be a closed 
subset of $(J_p \times \CC) \setminus J_2$ such that $f(T)\subseteq T$.  
Let $T_r$ denote the vertical neighborhood of $T$ of radius $r\geq 0$.  
Then there exists an $m_0 \in \NN$ and $r_0 >0$ such that for all 
$m \geq m_0$ and $r \leq r_0$, we have
$$ f^{m} (T_{r}) \subseteq T_{r/2}.  $$
\end{lem}

\medskip\noindent
When the conclusion of Lemma \ref{lem:trap} is satisfied, we say that 
$r_0$ is a {\bf trapping radius} (for $T$ under $f$).

\begin{proof}
Since $J_2$ is compact and $T$ is closed and disjoint from $J_2$, 
there is an $r_0>0$ with  $T_{r_0} \cap J_z = \emptyset$, for all $z \in J_p.$
Thus, Theorem \ref{thm:Cthm4.1} applied to each connected component
of the slice $T_{r_0} \cap ( \{ z \} \times \CC)$
yields that 
there exists an $m_0 \in \NN$ such that for all $m \geq m_0$ and $r \leq r_0$, 
we have  
 $$f^m( T_r ) \subseteq \Nvs (f^{m}(T), r/2) \subseteq T_{r/2}$$ 
since $f(T)\subseteq T$.  
\end{proof}

\begin{prop}  \label{prop:trapsets}
The sets 
  $$\Lambda, \Apt, \Acc, \A, PC_{J_p}$$
are disjoint from $J_2$ and have a trapping radius $r>0$.  
\end{prop}

\begin{proof}
First note that $f(\Lambda) = \Lambda$ (since $\Lambda$ is a union of basic sets) and $\Lambda\cap J_2 = \emptyset$, 
so Lemma \ref{lem:trap} clearly applies to $\Lambda$.  By Theorems~\ref{thm:Jthm3.1} 
and~\ref{thm:Jthm8.2},
we have $PC_{J_p} \cap J_2 = \emptyset$, and by definition, 
 $$\Apt \subseteq \Acc \subseteq \A \subseteq PC_{J_p}.$$
As accumulation sets always satisfy $f(A(X))\subseteq A(X)$, Lemma \ref{lem:trap}
applies to each of these postcritical sets.
\end{proof}


\medskip\noindent{\bf The unstable manifold $W^u(\Lambda)$ is an attractor.}
Recall the unstable manifold $W^u(\Lambda)$ is the set of all points $x\in \Ct$ for which there exists a backward orbit $x_{-k}$ converging to $\Lambda$.    

\begin{prop} \label{prop:Wattracts}
For any closed subset $X$ of 
$(J_p\times\CC)\setminus J_2$, its accumulation set $A(X)$ 
is contained in $W^u(\Lambda)\cap(J_p\times\CC)$.
\end{prop}

\noindent
The following proof is adapted from \cite{FShyp}, Propositions 4.2 and 4.5, 
where they prove an analogous result for ``s-hyperbolic" endomorphisms of $\Pt$.  

\begin{proof}
Let $E(J_2)$ be an open neighborhood of $J_2$ in $J_p\times\CC$ on which $f$ is expanding.
In particular, $f^{-1}(E(J_2))$ is strictly contained in $E(J_2)$.
Fix $x\not\in W^u(\Lambda)\cap (J_p\times\CC)$.  Its preimages 
$f^{-n}(x)$ must accumulate on the 
non-wandering set over $J_p$, $\Omega_{J_p} = \Omega \cap (J_p\times\CC)$. Since $\Omega_{J_p} = J_2 \sqcup \Lambda$, 
the existence of $E(J_2)$ implies
that the preimages $f^{-n}(x)$ will accumulate on $J_2$.  

Choose $N$ so that $f^{-N}(x) \subset 
E(J_2)$.  By continuity, there exists a neighborhood $V$ of $x$ in $J_p\times\CC$
so that $f^{-N}(V) \subset 
E(J_2)$.  It follows that $f^{-n}(V)$ converges uniformly to $J_2$ as $n\to\infty$.  

Fix a neighborhood $U$ of $W^u(\Lambda)\cap (J_p\times\CC)$ in $J_p\times\CC$ 
and a radius $R>0$.  Then there is an integer $N = N(U,R)$ so that 
$f^{-n} ((J_p\times D(0,R))\setminus U)$ is contained in $E(J_2)$ for all $n\geq N$.  
This in turn implies
that $f^n((J_p\times\CC)\setminus E(J_2))$  lies in $U\cup (J_p\times \{|w|\geq R\})$
for all $n\geq N$.  

Finally, let $X$ be any closed subset of $(J_p\times\CC)\setminus J_2$ and $A(X)$
its accumulation set.  Choose the expanding neighborhood $E(J_2)$ small enough 
so that $X\cap E(J_2) =\emptyset$.  
As $U$ and $R$ are arbitrary, we conclude
that $A(X)$ must be contained in $W^u(\Lambda)\cap(J_p\times\CC)$.
\end{proof}


As an immediate corollary to Proposition \ref{prop:Wattracts}, we obtain:  

\begin{cor}  \label{cor:AinW}
The accumulation set $\A$ satisfies $\A \subseteq W^u(\Lambda) \cap (J_p \times \CC)$.
\end{cor}

\medskip\noindent{\bf Points with bounded orbit.}
Let $K_{J_p}$ denote  $\bigcup_{z\in J_p} K_z = K \cap (J_p \times \CC)$. 

\begin{lem}  \label{lem:boundedWisL}
The unstable manifold $W^u(\Lambda)$ satisfies 
  $$W^u(\Lambda) \cap K_{J_p} = \Lambda.$$
\end{lem}

\begin{proof}
We clearly have the inclusion $\Lambda \subset W^u(\Lambda) \cap K_{J_p}$,
because $f(\Lambda)=\Lambda$ and $\Lambda\subset K_{J_p}$.  
Fix $x\in W^u(\Lambda) \cap K_{J_p}$. 
The the orbit of $x$ is bounded and its accumulation set $A(x)$ must 
lie in the non-wandering set of $f$.  Proposition 
\ref{prop:Wattracts} implies that $A(x)\subset W^u(\Lambda)$ so $A(x)$ is 
disjoint from $J_2$, but then this implies that $A(x)\subset \Lambda$.  Consequently
$x$ is in the stable manifold $W^s(\Lambda)$.  We have
$W^s(\Lambda)\cap W^u(\Lambda) = \Lambda$ by \cite[Proposition A.4]{MatSk},
so $x$ is in $\Lambda$.
\end{proof}

\begin{lem} \label{lem:WsKint}
$W^s(\Lambda) = K_{J_p}\setminus J_2$.
\end{lem}

\begin{proof}
The inclusion $W^s(\Lambda) \subset K_{J_p}\setminus J_2$ is clear from the 
definitions.  The converse follows from Proposition \ref{prop:Wattracts} and
Lemma \ref{lem:boundedWisL}.
\end{proof}

\begin{lem} \label{lem:AptisL}
We have $\Lambda = \Apt$.
\end{lem}

\begin{proof}
First note that the accumulation set $A(x)$ of any point $x \in J_p \times \CC$ will be a 
subset of $K_{J_p}$, or empty. 
Proposition \ref{prop:Wattracts} implies that $\Apt \subset W^u(\Lambda) \cap K_{J_p}$,
so Lemma \ref{lem:boundedWisL} gives  $\Apt \subset \Lambda$.

If $(z,w)$ is a period $n$ saddle point with $z\in J_p$, 
then $p^n(z) = z$.  The polynomial
$Q^n_{z}(w)$ is hyperbolic; hence its attracting periodic points, which are saddle points of $f$, attract critical points in 
$\{z\}\times\CC$.  Thus $(z,w) \in \Apt$.  
Hence $\Apt$ contains all saddle periodic points with $z \in J_p$.
As $\Apt$ is closed, it also contains the closure of all saddle periodic points 
over $J_p$, which is $\Lambda$.
\end{proof}


\medskip\noindent{\bf The immediate basin of a saddle set.}  
The following proposition should be interpreted as 
an analog of the one-dimensional result, 
where there is always a critical point in the immediate basin of an attracting cycle.
For any point $x\in (J_p\times\CC)\setminus J_2$, let $U_x$ denote the vertical 
Fatou component containing $x$; that is, the connected component of the 
vertical slice of $(J_p\times\CC)\setminus J_2$ containing $x$.

\begin{prop} \label{prop:Lbasin}
There exists $N\in\mathbb{N}$ 
so that $U_x$ contains a critical point of $f^N$ for all $x\in\Lambda$.
\end{prop}

We begin with a simple lemma about hyperbolic metrics on planar domains.
Let $D(a,r)$ denote a Euclidean disk with center $a$ and radius $r$.

\begin{lem}  \label{lem:hyperbolicmetric}
Suppose $U\subset\CC$ is a bounded domain with 
 $$D(x,r)\subseteq U \subseteq D(0,R).$$
Then the hyperbolic metric $\rho_U(z) |dz|$ of $U$ restricted to the disk $D(x,r/2)$
is comparable to the 
Euclidean metric.  Explicitly,
  $$\frac{2}{R} \leq \rho_U(z) \leq \frac{8}{3r}$$
for all $z\in D(x, r/2)$.  
\end{lem}

\begin{proof}
Recall that the hyperbolic metric on the unit disk $D(0,1)$ (with constant curvature -1)
is given by $2 |dz| / (1 - |z|^2)$.  The inclusions $D(x,r)\subseteq U\subseteq D(0,R)$ 
imply that 
  $$\frac{2R}{R^2 - z^2} \leq \rho_U(z) \leq \frac{2r}{r^2 - |z-x|^2}.$$
Restricting to $z\in D(x, r/2)$ gives the desired estimate.
\end{proof}

\medskip
\noindent{\bf Proof of Proposition \ref{prop:Lbasin}.}
For each $x\in \Lambda$, let $U_x$ be the vertical Fatou component containing $x$.
These domains are uniformly bounded, since all are contained in $K$.
Choose $R$ so that $U_x\subseteq D(0,R)$ for all $x\in \Lambda$.  Let $d_x$ denote the 
hyperbolic distance function on the domain $U_x$.  

Let $r(\Lambda)$ be a trapping radius for $\Lambda$, as guaranteed
by Proposition \ref{prop:trapsets}.  Then $\Nvs(x,r(\Lambda))\subseteq U_x$
for all $x\in\Lambda$, and $f^{mn}(\Lambda_r)\subseteq \Lambda_{r/2^n}$ for fixed $m=m(\Lambda)$, any $n\in\mathbb{N}$, and all $r\leq r(\Lambda)$.  Note that the spherical metric 
is comparable to the Euclidean metric on the bounded domain $D(0,R)$.  
In particular, for each $y\in U_x$ with 
$|x-y| = r(\Lambda)/2$, Lemma \ref{lem:hyperbolicmetric} implies that 
 $$d_{f^{mn}(x)} (f^{mn}(x), f^{mn}(y)) \leq \frac{8r(\Lambda)}{3r(\Lambda) 2^{n+1}} 
 	< \frac{2r(\Lambda)}{2R} \leq d_x(x,y)$$
when $n$ is sufficiently large.  Consequently, the proper holomorphic map 
$f^{mn}: U_x\to U_{f^{mn}(x)}$ is strictly contracting in the
hyperbolic metric for every $x\in\Lambda$.  It follows that there exists a critical point 
of $f^{mn}$ in $U_x$ for every $x\in\Lambda$.
\qed

\medskip
The above yields the following proposition, which will be useful for establishing $\A \supseteq W^u(\Lambda) \cap (J_p \times \CC)$.

\begin{prop} \label{prop:LinfNC}
For any $r>0$, there is an $N \geq 0$ such that
$$\Lambda \subset \bigcup_{n=0}^N \Nvs(f^n(C_{J_p}), r).$$
\end{prop}

\begin{proof}
Let $N_0$ be the integer given in Proposition \ref{prop:Lbasin}.  
Denote the critical locus of $f^{N_0}$ over $J_p$ by $C_{N_0} = \bigcup_{i=0}^{N_0-1} f^{-i}(C_{J_p})$.
The vertical expansion of $f^{N_0}$ implies that $C_{N_0}$ is uniformly bounded away from $J_2$.  As a critical set, note that $C_{N_0}$ contains no isolated points and is fiberwise continuous over $J_p$ (that is, $z\mapsto (C_{N_0})_z$ is continuous).  Recall that the saddle set $\Lambda$ is also uniformly bounded from $J_2$ and that $z\mapsto J_z$ is continuous over $J_p$.   

We can therefore construct a family of paths, one from each point $x$ in $\Lambda$ to a point in $C_{N_0}$ contained in the Fatou component $U_x$, such that all paths are uniformly bounded from $J_2$.  Indeed, suppose to the contrary that there is a sequence of points $x_n$ in $\Lambda$ so that any path joining $x_n$ to $C_{N_0}$ in $U_{x_n}$ has distance to the boundary less than $\delta_n>0$ with $\delta_n \to 0$.  Pass to a subsequence so that $x_n \to x$ in $\Lambda$.  Take any path $\gamma$ in $U_x$ joining $x$ to $C_{N_0}$.  Then $\gamma$ is some distance $\eps>0$ from the boundary of $U_x$ (and therefore from $J_2$).  By the continuity of $C_{N_0}$ in $z$ and the continuity of $z\mapsto J_z$, there are paths joining $x_n$ to $C_{N_0}$ in $U_{x_n}$ which remain distance $\eps/2$ from $J_2$ for all large $n$.  This contradicts the assumption.

Using Theorem~\ref{thm:Cthm4.1}, iterating forward by some $N>N_0$ iterates, the paths contract
uniformly, implying that for each $x\in\Lambda$, there exists a point of $\bigcup_{i=0}^N f^i(C_{J_p})$
in the Fatou component $U_x$ within spherical distance $r$ of $x$.  In other words, $\Lambda$
is contained in the vertical neighborhood about $\bigcup_{i=0}^N f^i(C_{J_p})$ of radius $r$.
\end{proof}


\medskip\noindent{\bf Equality of $\A$ and the unstable manifold.}  
Now we are ready to establish the final inclusion needed for 
the statement of Theorem \ref{thm:Achain}.  

\begin{lem} \label{lem:WinA}
We have $W^u(\Lambda) \cap (J_p \times \CC) \subseteq \A$.
\end{lem}

\begin{proof}
Let $r_0$ be a trapping radius for the postcritical set $PC_{J_p}$.  
Proposition \ref{prop:LinfNC} says that the set $U = 
\bigcup_{n=0}^N \Nvs( f^n(C_{J_p}), r_0)$ contains a vertical neighborhood
of $\Lambda$ of some small radius.  Since $z\mapsto C_z
= C_{J_p}\cap (\{z\}\times\CC)$ is continuous in the Hausdorff topology, 
continuity of $f$ guarantees that $U$ contains
a neighborhood $V$ 
of $\Lambda$ in the ambient space $J_p\times\CC$.  

Fix $(z_0, w_0) \in W^u(\Lambda) \cap (J_p \times \CC)$.  
Let $(z_{-k}, w_{-k})\in f^{-k}(z_0,w_0)$ be a prehistory tending to $\Lambda$. 
Then  $(z_{-k}, w_{-k})$ is in $V \subset U$ for all sufficiently large
$k$.  Let $(z_{-k}, y_k)\in \bigcup_{n=0}^N f^n(C_{J_p}) $ be the closest point to $(z_{-k}, w_{-k})$ 
in the fiber $\{z_{-k}\}\times\CC$.  
Because $r_0$ is a trapping radius, 
the image $f^k(z_{-k},y_k)$ is very close to $(z_0, w_0)$.  Letting $k\to\infty$,
we have $(z_0,w_0)\in\A$.  
\end{proof}

\medskip\noindent{\bf Proof of Theorem \ref{thm:Achain}.}
The inclusions $\Apt \subseteq \Acc \subseteq \A$ are clear from the definitions. 
Combine Lemma \ref{lem:AptisL} with Corollary \ref{cor:AinW} and Lemma
\ref{lem:WinA} to complete the proof.  
\qed

\section{Stability and perturbations of products}
\label{sec:fiberhomeo}

In this section, we prove Theorem~\ref{thm:SkewHM}, showing that for holomorphic families of Axiom A skew products, 
the holomorphic motion of $J_2$ must preserve the vertical fibration.  As an application, 
we see that if a polynomial skew product is 
in the same hyperbolic component as a product, then all fiber Julia sets are homeomorphic.  

\bigskip\noindent{\bf Holomorphic motions.}
Suppose $E$ is 
a subset of a complex manifold $X$.
As above, $D(a,r)$ denotes a disk in $\CC$ centered at $a$ with radius $r$. 
Then $\Phi: D(0,r) \times E \to X$ is a {\bf holomorphic motion} of $E$ if $\Phi$ is continuous and
\begin{enumerate}
\item	$\Phi(0,e) = e$ for all $e\in E$,
\item	$\Phi(\cdot, e)$ is holomorphic for each fixed $e\in E$,
\item 	$\Phi_\lam := \Phi(\lambda, \cdot)$ is injective for each fixed $\lambda\in D(0,r)$.
\end{enumerate}

\bigskip\noindent{\bf Uniform expansion on $J_2$.}
Let $\{f_\lam$: $\lam\in D(0,1)\}$ be a holomorphic family of polynomial skew products of the form 
 $$f_{\lam} (z,w) = (p_{\lam}(z), q_\lam(z,w));$$
that is, each $f_\lam$ is a polynomial skew product of $\Ct$ and the coefficients of 
$p_\lam$ and $q_\lam$ are holomorphic in $\lam$.  If 
$p_{\lam}$ is hyperbolic for each $\lam\in D(0,1)$, and if each $f_\lam$ is 
vertically expanding over $J_{p_\lam}$, we say that the family $\{f_\lam\}$
is {\bf uniformly expanding} on $J_2$.    
Note uniform expansion on $J_2$ is weaker than Axiom A.

\begin{thm} \cite[Theorem C]{MatHM} \label{thm:matHM}
Let $\{f_\lam: \lam\in D(0,1)\}$ be a holomorphic family of polynomial skew products which is 
uniformly expanding on $J_2$.  Then there exists an $r>0$ and a holomorphic 
motion $\Phi : D(0,r) \times J_2(f_0) \to \Ct$ such that 
such that $\Phi_\lam(J_2(f_0)) = J_2(f_\lambda)$ and 
$$f_{\lam} = \Phi_{\lam}\circ f_0\circ \Phi_\lam^{-1}$$ 
on $J_2(f_\lam)$ for all $\lam\in D(0,r)$.  
\end{thm}

We show the holomorphic motion $\Phi$ also preserves the vertical fibration:

\begin{thm} \label{thm:SkewHM}
Under the same hypothesis as in Theorem~\ref{thm:matHM} we also have
that $\Phi_{\lam}$ is a skew product
 $$\Phi_{\lam}(z,w) = (\varphi_{\lam}(z), \psi_{\lam}(z,w)),$$ 
where 
\begin{enumerate}
\item	$\varphi \colon D(0,r) \times J_{p_0} \to \CC$ is a holomorphic motion of $J_{p_0}$ such that 
		$\varphi_\lambda$ conjugates $p_0|J_{p_0}$ to $p_\lam|J_{p_\lam}$, and 
\item	for each $z\in J_{p_0}$, $\psi_\lam(z,w)$ defines a holomorphic motion 
		$\psi: D(0,r) \times J_z(f_0)\to \CC$
		such that $\psi_\lam(J_z(f_0)) = J_{\varphi_\lam(z)}(f_\lam)$.
\end{enumerate}
\end{thm}

\begin{proof}
Because the one-dimensional polynomials $\{p_\lam\}$ are hyperbolic, there exists a 
holomorphic motion $\varphi: D(0,r) \times J_{p_0}\to\CC$ conjugating
$p_0|J_{p_0}$ to $p_\lam|J_{p_\lam}$ by the one-dimensional theory (see e.g.\ \cite{McMullen:CDR}).
By the density of repelling cycles in $J_{p_0}$, the motion is uniquely determined.  

Let $(z_0,w_0)\in J_2(f_0)$ be a periodic point of $f_0$, so $(z_\lam, w_\lam) = 
\Phi_\lam(z_0,w_0)$ is a periodic point of $f_\lam$ of the same period.  Then 
the point $z_\lam$ must be a repelling periodic point of $p_\lam$ of 
constant period in $\lam$.  Consequently,
$z_\lam = \varphi_\lam(z_0)$ for all $\lam\in D(0,r)$.  
This holds for all periodic points $(z,w)$ of $f_0$ with $z=z_0$,
and by density of periodic points in the fiber Julia set 
$J_{z_0}(f_0)$, we obtain that the projection to
the first coordinate of $\Phi_\lam(z,w)$ is $z_\lam$ for all $w\in J_{z_0}(f_0)$.  Finally,
by the density of periodic points in $J_{p_0}$, continuity of $\Phi_\lam$, and the fact that $z \mapsto J_z$ is continuous over $J_p$, we obtain that 
$\Phi_{\lam}(z,w) = (\varphi_{\lam}(z), \psi_{\lam}(z,w))$ for some function $\psi_\lam$
and all $(z,w)\in J_2(f_0)$.  

Since $\Phi_\lam$ is a skew product, the motion therefore preserves the vertical fibration
of the skew products $f_\lam$.  The proof of (2) then follows immediately from the properties
of $\Phi$ as a holomorphic motion.
\end{proof}

\bigskip\noindent{\bf Perturbations of a product.}
Theorem~\ref{thm:SkewHM} applied to  a product $f_0(z,w) = (p(z), q(w))$ gives the following.

\begin{prop}  \label{prop:fiberJzhomeo}
Let $\{f_\lam: \lam\in D(0,1) \}$ be a holomorphic family of polynomial skew products which is 
uniformly expanding on $J_2$.
If $f_0(z,w) = (p(z), q(w))$ is a product, then for each $f_{\lam}$ with $\lam \in D(0,r)$, 
the fiber Julia sets $J_z(f_\lam)$ are homeomorphic to $J_q$ for all $z\in J_{p_\lam}$.  
\end{prop}

\begin{proof}
For the product $f_0$, all fiber Julia sets are equal to $J_q$.  
Let $\Phi: D(0,r) \times J_2(f_0)\to \CC^2$ be the holomorphic motion guaranteed 
by Theorem \ref{thm:matHM}.  The result follows immediately from property (2)
of Theorem~\ref{thm:SkewHM}.
\end{proof}

By Theorem \ref{thm:Jthm8.2}, a product $f_0(z,w) = (p(z), q(w))$ is Axiom A if and only if each of $p$ and $q$ are hyperbolic polynomial maps of $\CC$. Thus for products, uniform expansion of $J_2$ is equivalent to Axiom A.   

\begin{cor}  \label{cor:productcomponent}
Suppose $f(z,w) = (p(z), q(z,w))$ is an Axiom A polynomial skew product in the same hyperbolic component
as a product $f_0(z,w) = (p_0(z), q_0(w))$.  Then for all $z\in J_p$, the fiber Julia sets $J_z$ of $f$ are homeomorphic (to $J_{q_0}$).  
\end{cor}

\begin{proof}
Because $f$ and $f_0$ are in the same hyperbolic component, they can be connected 
by a chain of holomorphic motions, as guaranteed by Theorem \ref{thm:matHM}.  By
Proposition \ref{prop:fiberJzhomeo}, the fiber Julia sets are homeomorphic 
to $J_{q_0}$.  
\end{proof}

In the following section, the family of examples $F_a(z,w) = (z^2, w^2 + az)$ will show that 
the converse to Corollary \ref{cor:productcomponent} is false.  That is, the fiber Julia sets for these
maps are all homeomorphic, but for appropriate choices of 
the parameter $a$, they are not in the same hyperbolic component as a product.

\section{Axiom A skew products with $\Lambda = W^u(\Lambda)\cap (J_p\times \CC)$}
\label{sec:L=W}

In this section, we first provide some general conditions under which we have equality for all inclusions
listed in Theorem \ref{thm:Achain} (Theorem~\ref{thm:Lbigequiv}).  As a corollary, we see equality is preserved in hyperbolic components.  We also show that equality holds for Axiom A products and their perturbations (Proposition \ref{prop:products}).  Finally, we give an infinite family of distinct, non-product
Axiom A maps, for which equality holds:

\begin{thm}  \label{thm:Fa}
Let $F_a(z,w) = (z^2, w^2 + az)$ and $g_a(w) = w^2+a$, for each  $a\in\CC$. We have:
\begin{enumerate}
\item	$F_a$ is Axiom A if and only if $g_a$ is hyperbolic.
\item	If $F_a$ is Axiom A, then it is in the same hyperbolic component as 
		a product if and only if $g_a$ has an attracting fixed point.
\item If $F_a$ is Axiom A, then $\Lambda = W^u(\Lambda)\cap (J_p\times\CC)$.
\end{enumerate}
\end{thm}

\bigskip\noindent{\bf Criteria for $\Lambda = W^u(\Lambda) \cap (J_p \times \CC)$.}
In this subsection, let $f$ be an Axiom A polynomial skew product.
We now give necessary and sufficient conditions which guarantee equality for all inclusions in 
Theorem \ref{thm:Achain}.  
As usual, if $T \subset \Ct$ then $T_z = T \cap (\{z \} \times \CC)$, and we say $z \mapsto T_z$ is continuous if it is continuous in the Hausdorff topology.  Recall $D(a,r)$ is a Euclidean disk with center $a$ and radius $r$.

\begin{thm} \label{thm:Lbigequiv}
The following are equivalent:
\begin{enumerate}
\item[(a)]  $\Lambda = W^u(\Lambda) \cap (J_p \times \CC)$;
\item[(b)] $z \mapsto \Lambda_z$ is continuous for all $z \in J_p$;
\item[(c)] $z \mapsto K_z$ is continuous for all $z \in J_p$.
\end{enumerate}
\end{thm}

\noindent
The implication (b) $\Rightarrow$ (a) was inspired by Robinson's analogous statement for diffeomorphisms (\cite[Ch.~8, Theorem 6.2]{Rob}).

\begin{proof}
First we show that if $\Lambda=\emptyset$ then each of (a), (b), and (c) hold. Note (a) and (b) are satisfied vacuously.    Now, $\Lambda = \emptyset$ implies that $W^s(\Lambda) = \emptyset$, which
implies that $K_{J_p}= J_2$ by Lemma \ref{lem:WsKint}.  Since $J_2 = \bigcup_{z\in J_p} J_z$ (Theorem~\ref{thm:Jprop3.5}), this yields $K_z = J_z$ 
for all $z\in J_p$.  
Again applying Theorem \ref{thm:Jprop3.5}, $z \mapsto J_z$ is continuous over $J_p$ for $f$ Axiom A, we conclude that (c) holds.
Now assume that $\Lambda \neq \emptyset$.

\textbf{(a) $\Rightarrow$ (b):}
Suppose $z \mapsto \Lambda_z$ is not continuous over $J_p$.  
Since $\Lambda$ is closed, we must have upper semi-continuity.  
Let $(z_0,w_0) \in \Lambda$ be a point where lower semi-continuity fails.  
Then there is a sequence $z_n\to z_0$ in $J_p$ and a $\delta>0$ so that 
$\Lambda_{z_n} \cap D(w_0,\delta) = \emptyset$ for all $n$.  
The local stable manifold of the point $(z_0,w_0)$ lies in $K_{z_0}$ in the vertical
fiber, therefore the local unstable manifold (for any choice of prehistory)
of $(z_0,w_0)$ must be transverse to the fiber (see e.g.\ \cite[\S2]{FShyp}).
Thus $W^u(\Lambda)_z \cap D(w_0,\delta) \not= \emptyset$ for all $z$ near $z_0$.  
This implies that $W^u(\Lambda)\cap (J_p\times\CC) \not= \Lambda$.  

\textbf{(b) $\Rightarrow$ (c):}
Suppose $z \mapsto K_z$ is not continuous over $J_p$.   Again since $K$ is closed, we can 
assume that $z\mapsto K_z$ fails to be lower semi-continuous at a point $(z_0,w_0)\in K$.  
As $z\mapsto J_z$ is continuous, we must have  $w_0\in K_{z_0}\setminus J_{z_0}$,
and there exists a sequence $z_m\to z_0$ in $J_p$ so that $(z_m,w_0)\not\in K$.

By Lemma \ref{lem:WsKint}, $(z_0,w_0)\in W^s(\Lambda)$.  
Let $f^{n_k}(z_0,w_0)$ be a subsequence of iterates converging to a point $(a,b)\in \Lambda$  
as $k\to\infty$.  By continuity of $f$, for each fixed $k$ the images $(x_m^k, y_m^k) = 
f^{n_k}(z_m,w_0)$ converge to $f^{n_k}(z_0,w_0)$
as $m\to\infty$.  

By continuity of $z\mapsto J_z$, there is a neighborhood $U$ of $a$ in $J_p$ 
and a $\delta>0$ so that $U\times D(b,\delta)$ is disjoint from $J_2$.  
Choose $k$ large enough so that $f^{n_k}(z_0,w_0)$ lies $U\times D(b,\delta)$; then 
$(x_m^k, y_m^k)$ lies in $U\times D(b,\delta)$ for all $m$ sufficiently large.  But $(z_m,w_0)\not\in K$,
so the invariance of $K$ implies that $D(b,\delta)\cap K_{x_m^k} = \emptyset$.
By shrinking $U$, we obtain a sequence of points $x_m^k$ converging to $a$ such that  
$\Lambda_{x_m^k}\cap D(b,\delta)\subset K_{x_m^k}\cap D(b,\delta) = \emptyset$,
and we conclude that $z\mapsto \Lambda_z$ is not continuous.

\textbf{(c) $\Rightarrow$ (a):}
By continuity of $z\mapsto K_z$, and since $\Lambda\subset K\setminus J_2$ is closed, 
 there is a neighborhood 
$U$ of $\Lambda$ in $J_p\times\CC$ such that $U\subset K\setminus J_2$.  

Let $(z_0,w_0)$ be any point in $W^u(\Lambda) \cap (J_p \times \CC)$.  Then there is sequence
of preimages $(z_{-k}, w_{-k}) \in f^{-k}(z_0,w_0)$ with $(z_{-k},w_{-k})\in U$ for all sufficiently 
large $k$.  This implies that $w_{-k}\in K_{z_{-k}}\setminus J_{z_{-k}}$.  By complete invariance
of $K$, we obtain $w_0 \in K_{z_0}$.
Recalling that $W^u(\Lambda)\cap K_{J_p} = \Lambda$ by Lemma \ref{lem:boundedWisL}, 
we conclude that $(z_0,w_0)\in\Lambda$.
Therefore $W^u(\Lambda)\cap (J_p\times\CC) = \Lambda$.  
\end{proof}

\begin{cor}  \label{cor:fconnected}
If $J_z$ is connected for all $z \in J_p$, then $\Lambda = W^u(\Lambda)\cap(J_p\times\CC)$.
\end{cor}

\begin{proof}
\cite[Lemma 3.7]{MatSk} states that if $f$ is vertically expanding over $J_p$, and if for all $z \in J_p$, $J_z$ is connected, then $z \mapsto K_z$ is continuous.
\end{proof}


\medskip\noindent{\bf Stability under perturbations.}   
We are now ready to prove Proposition \ref{prop:products}, which states that Axiom A
products and their perturbations satisfy $\Lambda = W^u(\Lambda)\cap(J_p\times\CC)$.
We first give a general proposition about perturbations.

\begin{prop} \label{prop:Lcomponent}
Suppose $f_0$ and $f_1$ are in the same hyperbolic component.  Then
$\Lambda= W^u(\Lambda)\cap(J_p\times\CC)$ holds for $f_0$ if and only if 
it holds for $f_1$.  
\end{prop}

\begin{proof}
We apply the characterization of 
$\Lambda= W^u(\Lambda)\cap(J_p\times\CC)$ from Theorem \ref{thm:Lbigequiv}.
Suppose that $z\mapsto K_z$ is continuous for $f_0$.  By Theorems \ref{thm:matHM}
and \ref{thm:SkewHM}, there is a holomorphic motion of $J_2$ on a neighborhood of $f_0$ 
which restricts to a holomorphic motion of $J_z$ in every fiber over $J_{p_0}$.  By the Sullivan-Thurston $\lambda$-lemma (\cite{Sullivan:Thurston}), the motion of $J_z$ extends to a motion of the whole fiber $\{z\}\times\CC$, perhaps restricting the domain of the motion.  Consequently, the topology of the filled Julia set $K_z$ is unchanged for nearby maps $f_\lambda$ (i.e. the Fatou components move with $J_z$).  The continuity of $z\mapsto J_z$ for $f_\lambda$ then guarantees that we also have continuity of $z\mapsto K_z$ for $f_\lambda$.  

Note that if $z\mapsto K_z$ is discontinuous for $f_0$, then it must fail to be lower semi-continuous.  Combining this with the continuity of $z\mapsto J_z$, there exists a sequence $z_n\to z_0$ in $J_p$ and a component $U$ of the interior of $K_{z_0}$, so that for any fixed compact subset $V\subset U$, $V$ is disjoint from $K_{z_n}$ for all $n$ large.  The Sullivan-Thurston $\lambda$-lemma preserves this discontinuity under perturbation.  

Finally, we can connect $f_0$ to $f_1$ by a finite chain of holomorphic motions which preserve the equality $\Lambda= W^u(\Lambda)\cap(J_p\times\CC)$, proving the proposition. 
\end{proof}

\medskip\noindent{\bf Proof of Proposition \ref{prop:products}.}
Let $f(z,w) = (p(z),q(w))$ be an Axiom A product.  Then $K_z = K_q$ for 
all $z\in J_p$.  From Theorem  \ref{thm:Lbigequiv}, it follows that 
$\Lambda= W^u(\Lambda)\cap(J_p\times\CC)$ for $f$.  Then Proposition
\ref{prop:Lcomponent} shows that all maps in same hyperbolic component as 
$f$ have $\Lambda= W^u(\Lambda)\cap(J_p\times\CC)$.
\qed


\bigskip\noindent{\bf The family $F_a$.}  We now give the proof of Theorem \ref{thm:Fa}.  
Let  $F_a(z,w) = (z^2, w^2 + az)$, $g_a(w) = w^2 + a$, and $p(z) = z^2$.  
For each $x\in \CC$, set 
	$$S_x = \{(e^{2it}, x e^{it}): t\in [0,2\pi]\}.$$ 
From the definition of $F_a$, we have
\begin{equation} \label{Fa1}
F_a(e^{2it}, x e^{it}) = (e^{4it}, g_a(x) e^{2it}), \mbox{ so } F_a(S_x) = S_{g_a(x)}.
\end{equation}
Note that the critical locus of $F_a$ over $J_p$ is $C_{J_p} = S_0$.  Therefore,
\begin{equation} \label{Fa2}
PC_{J_p} = \overline{\bigcup_{n>0} S_{g_a^n(0)}}.
\end{equation}

We begin with a lemma on the structure of $J_2$.

\begin{lem} \label{lem:JFa}
For each $a\in\CC$, we have $J_2 = \bigcup_{x\in J_{g_a}} S_x$.
\end{lem}

\begin{figure}
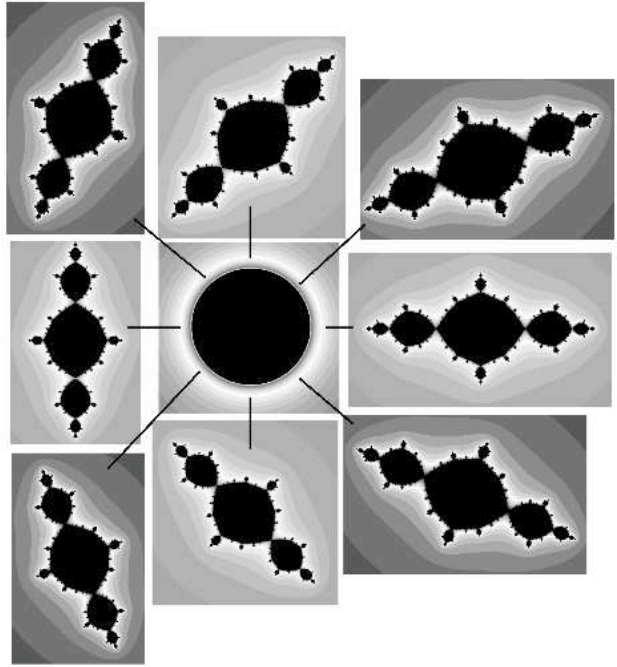

\begin{center}
  \drawfigtwbascl
\caption{\label{fig:twbas}
Let $F_a(z,w) = (z^2, w^2+az)$ be the family of Theorem~\ref{thm:Fa}. 
This figure depicts a collection of slices of the filled Julia set $K$, for  $a=-1$, with 
$K$ in black, and shades of gray distinguishing rate of escape level sets.  
In the center is the filled Julia set $K_p$ for the base polynomial $p(z)=z^2$. Circling around this are some fiber Julia sets, $K_z = K\cap (\{ z \} \times \CC)$, with lines drawn from the fiber $\{ z \} \times \CC$ to the corresponding point $z$ in $J_p$.  
Each $K_z$ is a rotation of the basillica (the Julia set of $w \mapsto w^2-1$), 
with rotation angle $\theta/2$ in fiber $z=e^{i\theta}$.
%
Figures~\ref{fig:twbas}, \ref{fig:aerocant}, and~\ref{fig:cantcircbas} were drawn with Fractalasm \cite{FAweb}.
}
\end{center}
\end{figure}

See Figure~\ref{fig:twbas}, which shows slices of $K$ for the map $F_{-1}$.


\begin{proof}
First observe that the vertical derivative of $F_a$ along the curve $S_x$ is 
given by
$$ q_{e^{2it}}'(xe^{it}) = 2x e^{it} = g_a'(x) e^{it}.  $$
If $x$ is a repelling periodic
point of $g_a$ of period $n$, then 
the vertical derivative of $F_a^n$ along the orbit
of $S_x$ satisfies $|(Q^n)'(w)| = |(g_a^n)'(x)| >1$.  Thus the iterates of $F_a^n$ cannot
be normal (in the fiber direction) 
along this orbit, so $S_x \subset J_2$, by Proposition~\ref{prop:Jprop2.1}.  Since $J_2$ is closed, $S_x \subset J_2$ for all $x\in J_{g_a}$.

Let $T=\bigcup_{x\in J_{g_a}} S_x$.  
For the reverse inclusion, note that $J_{z=1} = J_{g_a}$.  By invariance, we have
$J_z = T_z := T \cap (\{z\}\times\CC)$ for all $z$ along the backward orbit $p^{-n}(1)$.  These points are
dense in $S^1$.  By lower semi-continuity of $z\mapsto J_z$ (Proposition~\ref{prop:Jprop2.1}), we conclude that $J_z
\subset T_z$ for all $z$ in $J_p$.  Therefore $J_2 \subset T$.
\end{proof}

\begin{lem}  \label{lem:FaAxiomA}
$F_a$ is Axiom A if and only if $g_a$ is hyperbolic.
\end{lem}

\begin{proof}
Suppose $g_a$ is hyperbolic.  Then the postcritical set for $g_a$ remains a bounded
distance away from the Julia set $J_{g_a}$.  
The description of $PC_{J_p}$ for $F_a$ in Equation~(\ref{Fa2}) and Lemma \ref{lem:JFa}
imply that $PC_{J_p}$ is therefore disjoint from $J_2$. Consequently $F_a$ 
is vertically expanding over $J_p$.  
Note that the attracting set for the base map $p$ is $A_p = \{0\}$ and $q_0(w) = w^2$
for all $a\in\CC$.  
Therefore $F_a$ is also vertically expanding over $A_p$, so $F_a$ is Axiom A.

Conversely, if $g_a$ is not hyperbolic, then the postcritical set of $g_a$ is not
disjoint from $J_{g_a}$.  As $F_a| \{1\}\times\CC$ coincides with $g_a$, we conclude that 
$PC_{J_p}$ is
not disjoint from $J_2$.  Therefore $F$ is not vertically expanding over $J_p$ and
thus not Axiom A.
\end{proof}

\begin{lem}  \label{lem:Faproduct}
Suppose $F_a$ is Axiom A.  Then 
$F_a$ is in the same hyperbolic component as a product 
if and only if $g_a$ has an attracting fixed point.
\end{lem}

First, we recall some facts on the dynamics of the family $g_a(w) = w^2+a$ (see e.g., \cite[VIII.1]{Carleson:Gamelin}).
The {\bf Mandelbrot set} is 
	$$\mathcal{M} = \{ a \in \CC \colon J_{g_a} \mbox{ is connected} \}.$$  
For 
$a \in \CC \setminus \mathcal{M}$, the orbit of the critical point $0$ escapes to infinity, $g_a$ is hyperbolic, and $J_{g_a}$ is a Cantor set.  If $a$ is hyperbolic and in the interior of $\mathcal{M}$, then $g_a$ has an attracting cycle which attracts the orbit of the critical point $0$.

\begin{proof}
If $F_a$ lies in the same component as a product, then the product must be 
$H_a(z,w) = (z^2, w^2+a)$, because $F_a$ over the fixed point of $z^2$ 
in $J_p$ is given by the map $g_a(w) = w^2+a$ (and by Theorem \ref{thm:SkewHM}, the vertical fibration is preserved by holomorphic motions).  

First suppose that $g_a$ has an attracting fixed point. Then $g_a$ is in 
the same hyperbolic component as $g_0$ (the main cardiod of the
Mandelbrot set). Hence by Lemma~\ref{lem:FaAxiomA}, $F_a$ is in the 
same hyperbolic component as $F_0 = H_0$ which in turn is in the same 
hyperbolic component as $H_a$.  

Now suppose $g_a$ is hyperbolic with $g^n_a(0)\to\infty$.  
For the product we have
$J_2(H_a) = S^1 \times J_{g_a} $.  For the twisted map, Lemma \ref{lem:JFa}
implies that $J_2(F_a) = ([0,1] \times J_{g_a} )/
(0,x)\sim (1,-x)$, considered as a subset of $S^1 \times \CC$.  
If $F_a$ and $H_a$ are in the same hyperbolic component, then 
we can connect $F_a$ and $H_a$ by a chain of holomorphic motions which preserve
the vertical fibration (Theorem~\ref{thm:SkewHM}).  These 
motions therefore induce an isotopy from $J_2(H_a)$ to $J_2(F_a)$
within the ambient space $S^1 \times \CC$.  But this is impossible because connected
components of $J_2(H_a)$ are circles which project to $S^1$ with degree 1, while
connected components of $J_2(F_a)$ are circles which project to $S^1$ with degree
2, and these are in different homotopy classes.

Alternatively, suppose $g_a$ is hyperbolic with an attracting cycle of period $n>1$.  Suppose 
that $F_a$ and $H_a$ are in the same hyperbolic component.  Then by Theorem~\ref{thm:SkewHM}, $F_a$
and $H_a$ are conjugate on their Julia sets $J_2$ by a conjugacy which preserves
fibers over $J_p$.  In particular, the conjugacy must be the identity 
over $z=1$.  Let $\Gamma(H_a)$ (respectively $\Gamma(F_a)$)
be the subset of $J_2(H_a)$ (resp. $J_2(F_a)$) such that 
each slice $\Gamma(H_a)_z$ (resp. $\Gamma(F_a)_z$) is the boundary of the Fatou component in 
$\{z\}\times\CC$ containing the critical point $(z,0)$.  This subset of $J_2$ is dynamically
characterized in the following way: the slice $\Gamma(H_a)_z$ (resp. $\Gamma(F_a)_z$)
is the smallest connected subset of $J_z$ which is mapped to its
image with degree 2 by $H_a$ (resp. $F_a$).  Therefore, the conjugacy
must take $\Gamma(H_a)$ to $\Gamma(F_a)$.  Consequently, the conjugacy maps
the image $H_a(\Gamma(H_a))$ to $F_a(\Gamma(F_a))$.  Note however that 
the fibers of these image sets over $z=1$ do not coincide:  for the product we have a curve
winding around $w=a$ whereas for $F_a$ we have two curves winding around
$a$ and $-a$, contradicting the fact that the conjugacy is the identity over $z=1$.  
\end{proof}
 
\bigskip\noindent{\bf Proof of Theorem \ref{thm:Fa}.}
Combine Lemmas \ref{lem:FaAxiomA} and \ref{lem:Faproduct} for parts (1) and (2). 
When $a$ is not in the Mandelbrot set, the critical points escape so $\Lambda =
\Apt = \emptyset$ by Lemma \ref{lem:AptisL}
and we trivially obtain equality.   
For $a$ in a hyperbolic component in the Mandelbrot set, the map $F_a$ is 
connected, so part (3) for $f=F_a$ follows from Corollary \ref{cor:fconnected}.   
\qed

\bigskip
As a corollary, note by Proposition~\ref{prop:Lcomponent} that if $f$ is in the same hyperbolic component as any Axiom A $F_a$, then $\Lambda = W^u(\Lambda) \cap (J_p \times \CC)$ for $f$.


\begin{rem} 
By Theorem~\ref{thm:Jthm8.2}, since the extension of $F_a$ to the line at infinity is simply the map $\zeta \mapsto \zeta^2$, we have $F_a$ is Axiom A on $\Pt$ if and only if $F_a$ is Axiom A on $\Ct$.

Finally, suppose $F_a$ is Axiom A.  Note that if $g_a$ has disconnected Julia set, then $F_a$ satisfies (1) of Theorem~\ref{thm:examples},  while if $a$ lies in the Mandelbrot set, then $F_a$ satisfies (2) of Theorem~\ref{thm:examples}.
\end{rem}

\begin{rem} \label{rem:semiconj}
While the map $F_a$ does not lie in the same hyperbolic component of a product, it
should be noted that it is semiconjugate to the product $H_a(z,w) =
(z^2, w^2+a)$ via the map $\phi(z,w) = (z^2, zw)$.  That is, $F_a\circ \phi 
= \phi\circ H_a$, though this semiconjugacy does not extend regularly to $\Pt$.  
This was pointed out to the authors by Mattias Jonsson.  
\end{rem}

\medskip\noindent{\bf Sumi's example.}
Sumi has communicated to the authors the following very interesting example of a nonproduct Axiom A polynomial skew product satisfying $\Lambda = W^u(\Lambda) \cap (J_p \times \CC)$.  In his example, $J_p$ is a Cantor set, and all fiber $J_z$'s are connected. He constructs similar examples in \cite{Sumi}.

\begin{example} \label{exmp:sumi}
For any $R, \eps >0$ and  $n \in \NN$, let
$
p_R(z) = z^2 - R$, $p = p_{R}^n$, $h_{\eps}(w) = (w-\eps)^2 -1 + \eps$, and
  define $t_{n, \eps}(w)$ by
$h^n_{\eps}(w) = w^{2^n} + t_{n, \eps}(w)$.  

For appropriate choices of $\eps$ small, and $R, n$ large with $n$ even,
 the map
$$
f(z,w) = \(p(z), w^{2^n} + \(\frac{z+\sqrt{R}}{2\sqrt{R}}\) t_{n,\eps}(w) \)
$$
is an Axiom A polynomial skew product of $\Pt$ satisfying:
\begin{enumerate}
\item $C_{J_p} \subset K$;
\item $\Lambda = W^u(\Lambda) \cap (J_p \times \CC)$;
\item  $J_z$ is a Jordan curve, but not a quasicircle, for a.e.\ $z \in J_p$, in the maximal entropy measure of $p$;
\item $f$ is not in the same hyperbolic component as any product.
\end{enumerate}
\end{example}

\begin{proof}[Sketch of Proof.]
Axiom A and (1) can be proven by examining the postcritical set, in a similar way to our study of the family of maps found in Section~\ref{sec:matnonprodgen}. 

By (1) and Proposition~\ref{prop:Jprop2.3}, we have $J_z$ is connected for every $z \in J_p$, so by  Corollary~\ref{cor:fconnected} we get (2). 

$J_p$ is contained in two disks, $D = D(\sqrt{R}, r)$ and $-D$, for a small $r >0$ (such that $r/\sqrt{R} \to 0$ as $R \to \infty$).  The fiber maps $q_z$ for $z$ in $D$ are small perturbations of $h^n_{\eps}$, and the fiber maps $q_z$ for $z$ in $-D$ are small perturbations of $w \mapsto w^{2^n}$.  As a result, $J_{\beta}$ over the $\beta$-fixed point of $p$ (in $D$) is a quasi-basillica (not a Jordan curve), and $J_{\alpha}$ over the $\alpha$-fixed point of $p$ (in $-D$) is a quasi-circle.
Applying Lemmas 4.31 and 4.37 of \cite{Sumi} then yields (3).  

Finally, (4) follows from (3) and Corollary~\ref{cor:productcomponent}.
\end{proof}

\section{Axiom A  skew products with $\Apt \neq \Acc$}
\label{sec:aerocantor}

In this section, we show that the equality $\Apt = \Acc$ is preserved in hyperbolic components
(Proposition~\ref{prop:AptisAcc}),
and we construct an infinite family of Axiom A skew products in distinct hyperbolic components
for which $\Apt \neq \Acc$, yielding (3) of Theorem \ref{thm:examples}.

\medskip\noindent{\bf The $n$-airplane.}
Let $p_n(z) = z^2 + c_n$ be the unique quadratic polynomial with periodic critical point of least period $n$ and $c_n$ real.  
For example, $p_3(z) \approx z^2-1.75488$ is the ``airplane''. 
Then $\{ c_n \}$ is a sequence of real numbers descending to $-2$, 
and the Julia set $J_{p_n}$ is connected for each $n$.   
Let $\beta_n$ denote the $\beta$-fixed point of $p_n$, the point in $J_{p_n}$ with greatest real part.  

\begin{thm}  \label{thm:connectedbase}
Consider the sequence of skew products 
  $$f_n(z,w) = (p_n(z), w^2 + 2(2-z)).$$
For all sufficiently large $n$, $f_n$ is Axiom A and 
\begin{enumerate}
\item		$K\cap C_{J_{p_n}} = \{ (\beta_n, 0) \}$;
\item		 $\Lambda$  consists of a single fixed point in the 
		fiber $\{\beta_n\}\times\CC$;
\item		$J_z$ is disconnected for all $z\in J_{p_n}\setminus\{\beta_n\}$, while $J_{\beta_n}$
		is a quasicircle;
\item		$J_2$ is connected;
\item		$\Apt \not= \Acc$; 
\item	$f_n$ is not in the same hyperbolic component as a product; and
\item	$f_n$ is in the same hyperbolic component as $f_m$ if and only if $n=m$.
\end{enumerate}
\end{thm}

Each of the maps $f_n$ is a small perturbation of  $f_{\infty} (z,w)=(z^2-2, w^2+2(2-z))$.   Jonsson (\cite[Example 9.7]{MatSk}) shows $f_{\infty}$ is vertically expanding over $J_p$, and $f_{\infty}$ has the same connectivity properties as $f_n$ (i.e., (3) and (4) of the theorem), but $f_{\infty}$ is not Axiom A since the base is not hyperbolic.  Our examples $f_n$ are the first Axiom A examples
with such connectivity properties, which turn out to be the key to constructing an example with $\Apt \neq \Acc$. 
The outline of the proof of Theorem~\ref{thm:connectedbase} is based on \cite[Example 9.7]{MatSk}, but since $J_{p_n}$ is not contained in $\RR$, our case is of increased complexity.  

\begin{rem}
Because $f_n$ is a small perturbation of $f_\infty$, one might expect $f_n$
to be vertically expanding over $J_p$ because vertical expansion should be an 
open condition.  For general composition sequences, hyperbolicity is open in 
the $l^\infty$ topology on the space of sequences
\cite[Corollary 3.2]{MCrand}, but it is not open in the product topology, even with 
uniform coefficient bounds.  Consider the sequence of sequences $\{\{s_n^m: n\geq 1\}\}_{m\geq 1}$
given by
 $$s^m_n(w) = \left\{ \begin{array}{ll}  w^2 & \mbox{ for } n< m \\
							w^2 + 1/4 & \mbox{ for } n\geq m \end{array} \right.  $$
Then for each fixed $m$, the composition sequence $\{s^m_n\circ\cdots \circ s^m_1\}$ 
is not hyperbolic, because the critical point at $w=0$ is iterated towards the Julia set,
the locus of non-normality (see \cite[Theorem 1.3]{MCrand}).  On the other hand, the sequence converges in the 
product topology to $(w^2, w^2, w^2, \ldots)$ as $m\to \infty$, which is hyperbolic.  

Small perturbations of skew products correspond to small perturbations of fiberwise 
compositions in the product topology (with uniform bounds on the coefficients), not the 
$l^\infty$ topology (unless the nearby maps are conjugate).  
\end{rem}

\medskip

\begin{figure}
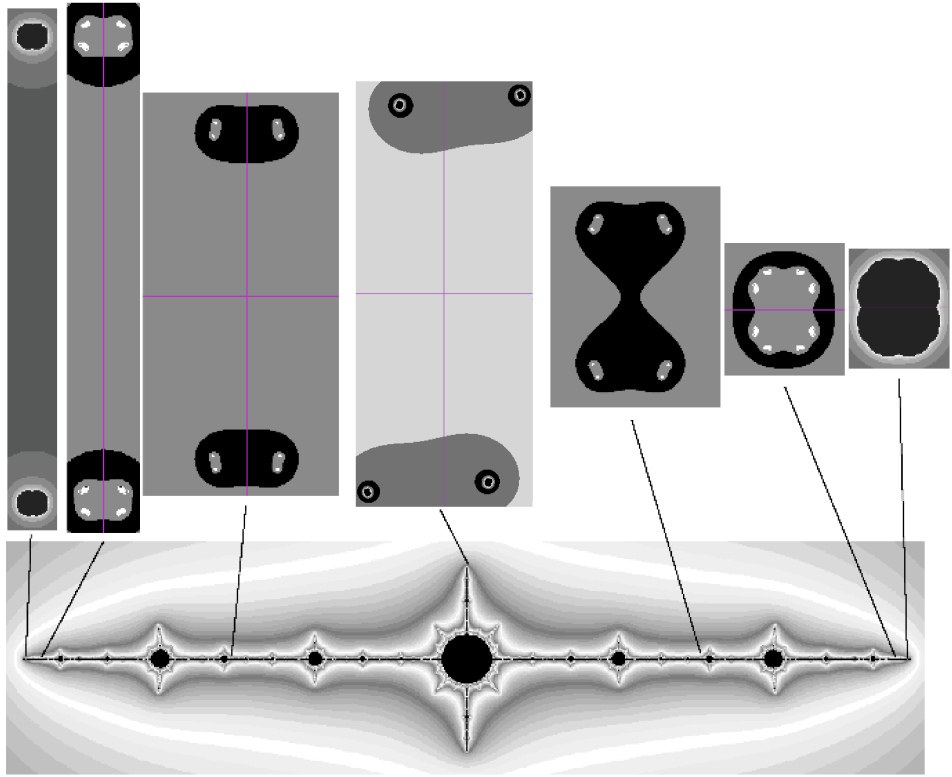

\begin{center}
  \drawfigaerocant
\end{center}
\caption{\label{fig:aerocant}
Let $f(z,w) =(p(z), w^2+2(2-z))$, for $p(z)=z^2-1.75488$.  Then $f(z) \approx f_3(z)$ from Theorem~\ref{thm:connectedbase}.  As in Figure~\ref{fig:twbas}, we show slices of $K$. The lower figure is $J_p$, and above (from right to left) are fibers: 
$z\approx 1.92, z=1.8, 1, z\approx 0.4i, z=-1, -1.8, z\approx -1.92$.  Note $K_{1.92}$ is a quasidisk, and $f$ maps $K_{-1.92}$ onto $K_{1.92}$ and $K_{0.4i}$ onto $K_{-1.92}$.
 }
\end{figure}

See Figure~\ref{fig:aerocant} for some slices of $K$ for a map of the type of 
Theorem~\ref{thm:connectedbase}. 

\medskip\noindent{\bf Maps with $\Apt = \Acc$.}

\begin{prop} \label{prop:AptisAcc}
Suppose $f_0$ and $f_1$ are in the same hyperbolic component.  Then the equality
$\Apt = \Acc$ holds for $f_0$ if and only if it holds for $f_1$. 
\end{prop}

We reduce the proof to the following lemma.

\begin{lem} \label{lem:AptAccchar}
If $f$ is Axiom A, then $\Apt = \Acc$ if and only if  
for any connected component $\mathcal{C}$ of $C_{J_p}$, we have either 
$\mathcal{C} \cap K = \emptyset$ or $\mathcal{C} \subset K$.
\end{lem}

\begin{proof}
Suppose $\Apt = \Acc$.  
Let $\mathcal{C}$ be a connected component of $C_{J_p}$.  
Since $\Lambda = \Apt$ by Lemma~\ref{lem:AptisL}, we also get $\Acc = \Lambda$.  
Thus $A(\mathcal{C})$ is either empty or contained in $\Lambda$.
If $A(\mathcal{C}) = \emptyset$, then $\mathcal{C} \cap K = \emptyset$. 
On the other hand, if $A(\mathcal{C}) \subseteq \Lambda$, then  
since $\Lambda \subset K$ and $\mathcal{C}$ is connected, by complete invariance of $K$ 
we must have $\mathcal{C} \subset K$.

Suppose for each connected component $\mathcal{C}$ of $C_{J_p}$,
we have either $\mathcal{C} \cap K = \emptyset$ or $\mathcal{C} \subset K$.
If $\mathcal{C} \cap K = \emptyset$ then $A(\mathcal{C}) = \emptyset$.  On  the other hand, if $\mathcal{C} \subset K$, then combining complete invariance of $K$ with 
Proposition~\ref{prop:Wattracts} and Lemma \ref{lem:boundedWisL}
yields $A(\mathcal{C}) \subseteq (W^u(\Lambda) \cap K_{J_p}) = \Lambda$.
Thus $\Acc \subseteq \Lambda = \Apt$, and the reverse inequality $\Apt\subseteq \Acc$
follows from the definitions.  
\end{proof}

\medskip\noindent{\bf Proof of Proposition \ref{prop:AptisAcc}.}
Let $f_0 = (p(z), q(z,w))$.  
Suppose $\{ f_{\lambda} = (p_\lam, q_\lam): \lambda \in D(0,1)\}$ is a holomorphic family of polynomial skew products which are uniformly expanding on $J_2$, with $\lambda=0$ giving the map $f_0$.  
Then by Theorems \ref{thm:matHM} and \ref{thm:SkewHM}, there is an $r>0$ and a holomorphic motion $\Phi \colon D(0,r) \times J_2(f_0) \to \Ct$ which conjugates the dynamics and  preserves the vertical fibration.  
In particular, $\Phi$ induces
one-dimensional holomorphic motions $\varphi$ of $J_p$ and $\psi_z$ 
of the fiber Julia sets $J_z$ for each $z\in J_p$.

For $\lam$ sufficiently small and each $z\in J_p$, the critical points of $(q_\lam)_{\varphi_\lam(z)}$
are close to the critical points of $q_z$.  If a critical point of $q_z$ lies
in $K_z\setminus J_z$, then the holomorphic motion $\psi_z$ of $J_z$ forces
nearby points to lie in $K_{\varphi_\lam(z)}\setminus J_{\varphi_\lam(z)}$ for all nearby $\lam$.  
Similarly, a critical point in $\CC\setminus K_z$ must remain 
in $\CC\setminus K_{\varphi_\lam(z)}$ under perturbation.  

Furthermore, the motion $\varphi$ of $J_p$ ensures that connected components of $C_{J_p}$
are uniformly close to connected components of $C_{J_{p_\lam}}$ for all 
$\lam$ sufficiently small.  
Let $\mathcal{C}$ be a connected component of $C_{J_p}$.  By Lemma \ref{lem:AptAccchar},
either $\mathcal{C}\cap K = \emptyset$ or $\mathcal{C}\subset K$.  Therefore, we have 
$\mathcal{C}\cap K = \emptyset$ or $\mathcal{C}\subset K$ for all connected components 
of $C_{J_{p_\lam}}$ for the maps $f_\lam$.  

Finally, connect $f_1$ to $f_0$ by a closed path in the hyperbolic component.  This path can be covered by a finite collection of overlapping open sets, on which the relation of $K$ to connected components $\mathcal{C}$ of $C_{J_p}$ as described above is constant.  The main result then follows from Lemma \ref{lem:AptAccchar}.  
\qed

\medskip
The remainder of the section is devoted to the Proof of Theorem~\ref{thm:connectedbase}.

\bigskip\noindent{\bf Proof of Theorem~\ref{thm:connectedbase}.}
Let $f_n(z) = (p_n(z), w^2 + 2(2-z))$ where $p_n$ is the $n$-airplane defined at the beginning of the section. 
Our most difficult task in this proof will be to establish that $f_n$ is Axiom A for $n$ sufficiently large.  
Since $p_n$ is a hyperbolic polynomial, by Theorem~\ref{thm:Jthm8.2}, we need only show 
$f_n$ is vertically expanding over $A_{p_n}$ and $J_{p_n}$.   To check vertical expansion we will
apply Theorem~\ref{thm:Jthm3.1}, and show the postcritical set over $A_{p_n} \sqcup J_{p_n}$ is disjoint from the fiber Julia sets.

Hence, our first step is to provide two lemmas, locating first (in Lemma~\ref{lem:Knconverges}) the base filled Julia set $K_{p_n}$ (which contains both $J_{p_n}$ and $A_{p_n}$), and then (in Lemma~\ref{lem:claim0}) the fiber $J_z$'s for $z \in K_{p_n}$.
As usual, we let $D(w,r)$ denote the open disk
in $\CC$ centered at $w$ with radius $r$.  
 For $w\in \CC$, we denote a closed rectangle
around $w$ by 
$$S(w, r_1, r_2) = 
 \{ x+iy \colon \abs{x- \Re(w)} \leq r_1, \abs{y-\Im(w)} \leq r_2 \},$$
and let $B(w, r) = S(w, r, r)$.  

 
\begin{lem} \label{lem:Knconverges}
There is a sequence $\eps_n \to 0$ such that 
$K_{p_n} \subset [-2, 2] \times [-\eps_n, \eps_n]$.
\end{lem}

\begin{proof}
Let $G_c(z) = \lim_{n\to\infty} 2^{-n} \log^+|p_c^n(z)|$ be the escape-rate function 
of $p_c(z) = z^2+c$.  Then $G_c$ is continuous as a function of both $z$ and $c$ (see e.g.\ 
\cite[VIII: Theorem 3.3]{Carleson:Gamelin}.  As a consequence, 
the mapping $c \mapsto J_c$ is lower semi-continuous (because $J_c =  \operatorname{supp} \Delta G_c$), 
while $c \mapsto K_c$ is upper semi-continuous (because $K_c = \{G_c=0\}$).
Hence for $c=-2$, we have $[-2,2]  = J_{-2} \subseteq  \lim_{c \to -2} J_c 
\subseteq  \lim_{c \to -2} K_c \subseteq K_{-2} = [-2,2]$.
Thus $c \mapsto K_c$ is continuous at $c=-2$.

Note also that the logarithmic capacity of $K_c$ is 1 for all $c$ (see e.g.\ \cite[VIII: Theorem 3.1]{Carleson:Gamelin}).
Since $c_n$ is in the Mandelbrot set, $K_{p_n} \subset \overline{D(0,2)}$, because any connected
set of logarithmic capacity 1 has diameter bounded by 4, and $K_{p_n}$ is symmetric about the origin.
Thus there is a sequence $\eps_n \to 0$ such that $K_{p_n} \subset [-2, 2] \times [-\eps_n, \eps_n]$.
\end{proof}


Recall that the fiber map $q_z(w) =
w^2 + 2(2-z)$ is independent of $n$.

\begin{lem}  \label{lem:claim0}
There is an $n_0$ such that for all $n \geq n_0,$ and all $z \in K_{p_n}$, we have
\begin{enumerate}
\item[(i)]  $q_z(\CC\setminus D(0,3.5)) \subset \CC\setminus D(0,3.75)$,
and 
\item[(ii)] $K_z \subset D(0,3.5)$.  
\end{enumerate}
\end{lem}

\begin{proof}
By Lemma~\ref{lem:Knconverges}, there is an $n_0$ such that $\eps_n \leq 1/4 $ 
if $n \geq n_0$.   
Then for $z \in K_{p_n}$ and $|w|\geq 3.5$, we have 
$|q_z(w)| = |w^2 + 2(2-z)| \geq |w|^2 - 2|2-z| \geq |w|(|w| - 2(4+\eps_n)/|w|) \geq
|w| (7/2 - 17/7) = 15|w|/14 \geq 15*3.5/14=3.75$, proving (i).  In particular, $q_z(w) \geq \frac{15}{14}|w|$, so the point $(z,w)$ escapes to infinity, proving (ii). 
\end{proof}


Now we can easily show vertical expansion over $A_{p_n}$.

\begin{lem}  \label{lem:claim1}
$f_n$ is vertically expanding over the attracting cycle $A_{p_n}$ of $p_n$.  
\end{lem}

\begin{proof}
Note that $A_{p_n}$ is real, and note $C_{A_{p_n}} = A_{p_n} \times \{ 0 \}$.  Let $x \in A_{p_n}$.  We show $(x, 0)$ escapes.  Let $(x_k, y_k) = f^k(x, 0)$.   Then for some $j \in \{ 0, \ldots, n \}$, we have $x_j = p_n^j(x) = c_n <0$.    
But then $y_{j+1} = y_j^2 + 2(2-x_j) \geq 4$, since $y_j$ is real.  Hence 
by Lemma~\ref{lem:claim0}, 
$y_{j+1} \notin K_{x_{j+1}}$, so $y_k \notin K_{x_k}$ for all $k$.  Thus $f$ is vertically expanding over $A_{p_n}$.
\end{proof}

Establishing vertical expansion over $J_{p_n}$ is the work of the next couple of pages.  Here is an overview. Note $C_{J_{p_n}} = J_{p_n} \times \{ 0 \}$.  We follow the outline of \cite[Example 9.7]{MatSk}.  
Fix a small $r \in (0, 1/16]$, setting $r=1/16$ suffices for our proof.  
The idea is that first, if $z \in J_{p_n}$ and $\Re(z) \leq 0$, then a small neighborhood of the real axis in the fiber $\{z \} \times \CC$ immediately escapes $D(0, 3.5)$ (which contains $K_z$ by Lemma~\ref{lem:claim0}).  
Next, if $z \in J_{p_n} \setminus \{ \beta_n \}$ and $\Re (z) > 0$, then the orbit $z$ under $p_n$ marches into $\{z \colon \Re(z) < 0\}$; further, for fixed $r>0$, there is a uniformity in the number of iterates it takes for the orbit of any $z \in J_{p_n} \setminus B(2,r)$ to land in $\{ \Re(z) < 0 \}$, for all sufficiently large $n$ (Lemma~\ref{lem:claim2}).  Combining the previous two ideas yields that critical orbits over $J_{p_n} \setminus B(2,r)$ escape and uniformly avoid $K$ (Lemma~\ref{lem:claim3}).  
Then in fibers over $J_{p_n} \cap B(2,r)$, we show a small neighborhood of the origin is mapped into itself (Lemma~\ref{lem:claim4}).  Hence critical points over $z$ remain near the origin (and in $\CC \setminus J_z$) as long as $z$ remains in $B(2,r)$, then once $z$ lands outside of $B(2,r)$, the previous case shows the orbits escape (Lemma~\ref{lem:claim5} (i)).  We add to this the dynamics of $w^2$ in the $\beta_n$-fiber (Lemma~\ref{lem:claim5} (ii)), to show the critical orbits in $J_{p_n} \cap B(2,r)$ avoid $J_2$ (Lemma~\ref{lem:claim5} (iii)).
Finally, we combine the above to get critical orbits over $J_{p_n}$ are uniformly bounded away from $J_2$ for sufficiently large $n$ (Lemma~\ref{lem:claim6}).

The following two lemmas locate the postcritical set in $J_{p_n} \setminus B(2,r)$. We first make a statement about the dynamics in the base (Lemma~\ref{lem:claim2}), then apply to the fibers (Lemma~\ref{lem:claim3}).  Let $n_0$ be chosen as in Lemma \ref{lem:claim0}.

\begin{lem} \label{lem:claim2}
There exist $N \geq 1$ and $n_1 \geq n_0$ such for all $n \geq n_1$ and $z \in 
{J_{p_n} \setminus B(2, r)}$, there is a $0\leq j < N$ with $\Re(p_n^j(z)) \leq 0$.
\end{lem}

\begin{proof}

Let $e_n : S^1 = \RR/\ZZ \to J_{p_n}$ be the external ray landing map for $p_n$.  
Since $p_n$ is hyperbolic, $e_n$ is well-defined and continuous, and it 
semi-conjugates $p_n$ on $J_{p_n}$ with angle doubling on $S^1$.  Recall that $\beta_n =e_n(0)$
is the $\beta$-fixed point of $p_n$.

First, we show that
for any $r  > 0$, there is a $\theta>0$ and $n_1 \geq n_0$ so that 
\begin{equation} \label{eqn:finiteN}
(-\theta, \theta) \subset e_n^{-1}(J_{p_n}\cap B(2,r))
\end{equation}
for all $n \geq n_1.$ 

As in the proof of Lemma \ref{lem:Knconverges}, the continuity of the 
escape-rate function $G_c(z)$ in both $c$ and $z$ implies that 
the harmonic measure $\mu_c = (2\pi)^{-1} \Delta G_c$ for the 
filled Julia set $K_c$ of  $z^2 +c$ is weakly
continuous in $c$.  Furthermore, the measure 
$\mu_c$ coincides with the push-forward of 
Lebesgue measure on $S^1$ by the external ray landing map.

Fix $r>0$ and consider $B = B(2,r) = [2-r, 2+r]\times [-r,r]$. 
Let $g : B \to [0,1]$ be a bump function supported in $B$ with $g= 1$ on the box
$B(2, r/2)$.  Then 
$$\mu_{-2}(B) > \int_B g \,d\mu_{-2} > C$$
for some $0< C< 1$.  This implies that at least $C$ of the total angle
lands in $B$.  By continuity of the landing map $e_\infty: S^1\to J_{-2}= [-2,2]$ 
for $z^2-2$ and the $z\mapsto\bar{z}$ symmetry of $J_{-2}$, we have that 
$e_{\infty}^{-1}(B)$ contains the interval $(-C/2, C/2)$. 

By weak continuity, the integral $\int_B g \,d\mu_c$ varies continuously in $c$.  Therefore 
$$\mu_{c_n}(B) \geq \int_B g \,d\mu_{c_n} > C$$
for all sufficiently large $n$.  
For these $n$, the total angle landing in $B$ is bounded
below by $C$, and because $c_n$ is real, we have maintained the 
$z\mapsto \bar{z}$ symmetry of $J_{p_n}$.  Therefore, 
$e_n^{-1}(B)$ contains the interval $(-C/2, C/2)$ for all $n$ sufficiently
large.  
Setting $\theta = C/2$ yields (\ref{eqn:finiteN}).

Finally, because the landing map $e_n$ defines a semiconjugacy between angle doubling on $S^1$ 
and $p_n$ on $J_{p_n}$, there exists a uniform $N$ so that the finite orbit 
$\{z, z_1 = p_n(z), \dots, z_{N-1} = p_n^{N-1}(z)\}$ contains an element with $\Re(z_j)\leq 0$
for all $z\in J_{p_n}\setminus B(2,r)$ and all $n\geq n_1$. 
\end{proof}


For given skew product $f_n$, 
let $Q_{n,z}^N$ denote the composition 
$q_{p_n^{N-1}(z)}\circ \cdots q_{p_n(z)}\circ q_z$ of fiber maps (recall the fiber map $q_z$ is
independent of $n$).  
Recall by Lemma~\ref{lem:claim0}, we know $K_z \subset D(0, 3.5)$ for $z \in K_{p_n}$.

\begin{lem}  \label{lem:claim3}
Let $N$ be given by Lemma \ref{lem:claim2}.  
There exist $n_2\geq n_1$ and $\delta >0$  so that 
 $$Q_{n,z}^N\left(\{ w: \abs{\Im(w)} \leq \delta \} \right) \cap D(0,3.5) = \emptyset$$
for all $z \in J_{p_n} \setminus B(2,r)$ and $n\geq n_2$.  
\end{lem}

\begin{proof}
By Lemma \ref{lem:Knconverges}, there exists $n_2\geq n_1$ so that 
	$$\eps_n < \frac{1}{N8^N}$$
for all $n\geq n_2$.   Choose $\delta$ with $0<\delta< 1/8^N$.  

Fix $n\geq n_2$ and $z\in J_{p_n}\setminus B(2,r)$.  Let $z_k = p_n^k(z)$ denote
the orbit of $z$.  Fix $w\in\CC$ with $|\Im w|\leq \delta$ and let $(z_k,w_k) = f_n^k(z,w)$.

Let $j\in \{0, \ldots, N-1\}$ be the least integer such that $\Re z_j\leq 0$.  
If $|\Re w_k|\geq 4$ for some $k \in \{ 1, \ldots , j\}$, then we conclude by Lemma
\ref{lem:claim0} that 
$w_n \notin D(0, 3.5)$ for all $n \geq k$, hence $w_N \notin D(0, 3.5)$.
Thus, we may assume
that $|\Re w_k| < 4$ for all $k\leq j$.  

From the formula for $f_n$, we have that $\Im w_1 = 2\Re w \Im w - 2 \Im z$, so that
$|\Im w_1| \leq 8 |\Im w| + 2\eps_n$.  By induction we obtain
	$$|\Im w_j| \leq 8^j |\Im w| + 2\eps_n \sum_{k=0}^{j-1} 8^k < 8^j \delta + 2j8^j \eps_n < \frac12$$
with the final inequality by our choices of $\delta$ and $n$.  Thus, 
	$$\Re w_{j+1} = (\Re w_j)^2 - (\Im w_j)^2 + 4 - 2\Re z_j \geq 4 - (\Im w_j)^2 > 3.5.$$
We conclude that $w_{j+1}\not\in D(0,3.5)$, so 
$w_N\not\in D(0,3.5)$ by Lemma \ref{lem:claim0}.  
\end{proof}

Next, we analyze the postcritical set over $J_{p_n} \cap B(2,r)$.  

\begin{lem}  \label{lem:claim4}
For any $\delta' < 1/4$, there exists $n_3$ so that for all $n\geq n_3$ and $z\in J_{p_n}\cap B(2,r)$,
$q_z( S(0, 1/4, \delta'))$ is contained in the interior of $S(0, 1/4, \delta')$.  
\end{lem}

\begin{proof}
Fix $\delta' < 1/4$.  By Lemma \ref{lem:Knconverges}, we 
may choose $n_3$ so that $\eps_n< \delta'/4$ for all $n\geq n_3$.

Fix  $z\in J_{p_n}\cap B(2,r)$, $w\in S(0, 1/4, \delta')$, and $n\geq n_3$, 
and let $(z_1, w_1) = f_n(z,w)$.  Then 
  $$|\Im w_1| = |2 \Re w \Im w - 2 \Im z| \leq \frac12 \delta' + 2\eps_n < \delta',$$
and 
  $$|\Re w_1| = |(\Re w)^2 - (\Im w)^2 + 4 - 2 \Re z| \leq \frac{1}{16} + (\delta')^2 + 2r < \frac{1}{4},$$
because $r = 1/16$.  
\end{proof}

Let $\delta$ and $n_2$ be given by Lemma \ref{lem:claim3}, and let $n_3$ be given 
by Lemma \ref{lem:claim4}.  We may assume $\delta<1/4$ and $n_3\geq n_2$.

\begin{lem}  \label{lem:claim5}
For all $n \geq n_3$, we have 
\begin{itemize}
\item[(i)]  $S(0, 1/4, \delta) \cap K_z = \emptyset,$ for all $z \in J_{p_n}  \setminus \{ \beta_n \}$, and 
\item[(ii)] $S(0,1/4, \delta) \cap J_{\beta_n} = \emptyset$; hence
\item[(iii)]
$J_{p_n}\times S(0,1/4,\delta) \subset (J_{p_n} \times\CC)\setminus J_2.$
\end{itemize}  

\end{lem}

\begin{proof} 
First note that any point $z\in (J_{p_n}\cap B(2,r))\setminus\{\beta_n\}$ leaves $B(2,r)$ 
after some number of iterates of $p_n$.  Therefore,  
 (i) follows immediately from Lemmas \ref{lem:claim4}, 
\ref{lem:claim3}, and \ref{lem:claim0}.  Statement (ii) follows because 
$p(\beta_n)=\beta_n$, and Lemma \ref{lem:claim4} shows that the iterates 
of $q_{\beta_n}$ form a normal family on $S(0,1/4, \delta)$, hence $S(0, 1/4, \delta) \subset  \CC \setminus J_{\beta_n}$ by Proposition~\ref{prop:Jprop2.1}. 
Recalling that $J_2 = \overline{ \bigcup_{z\in J_p} J_z}$, (i) and (ii) yield (iii).
\end{proof}

Finally, we show how to combine the above to show the critical orbits over $J_{p_n}$ avoid $J_2$, giving us vertical expansion over $J_{p_n}$.

\begin{lem}  \label{lem:claim6}
For all $n\geq n_3$, $f_n$ is vertically expanding over $J_{p_n}$. 
\end{lem}

\begin{proof}
We analyze the postcritical set of $f_n$.  
For $z\in J_{p_n}$, let $(z_k, w_k) = f_n^k(z,0)$.  

For the case $z = \beta_n$, the orbit $(z_k,w_k)$ 
lies in $S(0,1/4,\delta)$ for all $k$, by Lemma \ref{lem:claim4}.  
By Lemma~\ref{lem:claim5} (iii), the orbit $(z_k,w_k)$ is uniformly bounded away from $J_2$.  

For $z\in (J_{p_n}\cap B(2,r))\setminus\{\beta_n\}$, let $m$ be the least integer such that $z_m\not\in B(2,r)$.
By Lemma \ref{lem:claim4}, we have $w_k$ in the interior of $S(0, 1/4, \delta)$ for all $k\leq m$,
so by Lemma~\ref{lem:claim5} (iii), $w_k$ lies a definite distance away from $J_2$.  
Complete invariance of $K$ (and $J_2$) and the uniform $N$ 
in the statement of Lemma \ref{lem:claim3} implies that all $(z_k, w_k)$ for $k > m$ are uniformly bounded 
away from $K$ (and $J_2$), since $J_2 \subset K \subset J_{p_n}\times D(0,3.5)$ (Lemma~\ref{lem:claim0}).   
The previous sentence also applies to critical points $(z,0)$ when $z\in J_{p_n}\setminus B(2,r)$.  
\end{proof}

That $f_n$ is Axiom A for $n$ sufficiently large follows from the fact that $p_n$ is hyperbolic, Lemmas \ref{lem:claim1} and 
\ref{lem:claim6}, and Theorems~\ref{thm:Jthm3.1} and~\ref{thm:Jthm8.2}.

Finally, we turn to statements (1)--(7) of the theorem.  Lemma \ref{lem:claim5} shows that all critical points except the 
one in the fiber over $\beta_n$ escape, and it is clear from the construction that $(\beta_n,0)$
is attracted to a fixed point in its fiber.  Recalling that $\Lambda = \Apt$ from Lemma \ref{lem:AptisL},
we find that $\Lambda$ is precisely this fixed point. 

Because the critical points $(z,0)$ escape for $z\not=\beta_n$, the Julia sets $J_z$ are 
disconnected for $z\in J_{p_n}\setminus\{\beta_n\}$ (Proposition~\ref{prop:Jprop2.3}),
while $J_{\beta_n}$ is a quasicircle because $q_{\beta_n}(w)$
is a small perturbation of $w^2$.  As the base $J_{p_n}$ is connected,  it follows that
$J_2(f_n)$ is connected \cite[Lemma 6.7]{MatSk}.  

Next, since $J_{p_n}$ is connected, $C_{J_{p_n}} = J_{p_n} \times \{ 0 \}$ is a single connected component.  Hence 
$\Acc = \A$.   But we showed the critical point $(\beta_n, 0)$ is bounded, while the rest escape.  
Hence Lemma~\ref{lem:AptAccchar} yields $\Apt \neq \Acc$.

Proposition \ref{prop:AptisAcc} says the equality $\Apt=\Acc$ is preserved in hyperbolic
components, and we know equality holds for products (Proposition \ref{prop:products}); therefore,
$f_n$ is not in the same hyperbolic component as a product.  

Finally, $f_n$ and $f_m$ are in distinct hyperbolic components for $n\not=m$, because $p_n$ and $p_m$ are in distinct hyperbolic components, and holomorphic motions
of $J_2$ induce motions of the base by Theorem \ref{thm:SkewHM}.  
\qed

\begin{rem}
A map $f_n$ from Theorem~\ref{thm:connectedbase} which is Axiom A in $\Ct$ extends 
to the line at infinity as the map $\zeta \mapsto \zeta^2$, hence $f_n$ is also Axiom A on $\Pt$.  
Such a map satisfies (3) of Theorem~\ref{thm:examples}.
\end{rem}

\section{Axiom A skew products with $\Acc \neq \A$}
\label{sec:matnonprodgen}

In this section, we construct an infinite family of Axiom A
skew products in distinct hyperbolic components
which satisfy $\Acc \neq \A$, giving (4) of Theorem~\ref{thm:examples}.

\begin{thm} \label{thm:s1s2}
Given any two hyperbolic, monic polynomials $s_1, s_2:\CC\to\CC$ of degree $d\geq 2$, 
and positive integers $k_1 + k_2 = d$, there exists an Axiom A polynomial skew product $f(z) = (p(z), q(z,w))$ 
such that 
\begin{enumerate}
\item		$J_p$ is a Cantor set, with two disjoint, forward-invariant compact subsets 
		$L_1$ and $L_2$ such that $p|{L_i}$ is conjugate to the one-sided full shift on $k_i$
		symbols;
\item		$C_{J_p} \cap K \subset (L_1 \cup L_2) \times \CC$, 
		thus $J_z$ is disconnected if $z \in J_p \setminus (L_1 \cup L_2) $;
\item		for each $i = 1,2$, the restriction $f|(L_i\times\CC)$ is a small perturbation of the 
		product $(p|L_i) \times s_i$; 
\item		$\Lambda = \Lambda_1 \cup \Lambda_2$, where $\Lambda_i \subset (L_i \times \CC),$
		and $\Lambda_i \neq \emptyset$ if and only if not all critical points of $s_i$ escape; and
\item		if $\Lambda_i \neq \emptyset$ for either $i=1$ or $2$, then $\Acc \neq \A$ and
 the hyperbolic component containing $f$ does not contain a product.
\end{enumerate}

\end{thm}

Our construction is inspired by Proposition 3.8 and Example 3.9 of \cite{MatDil}, 
where the authors provide examples of polynomial skew products which are Axiom A on $\Pt$ 
and have ``nonterminal" (not minimal) basic saddle sets (so therefore are not in the same hyperbolic component as any product).
Their fiber maps are derived from a combination of $s_1(w)=w^d$ and $s_2(w) = w^d + R$ for a large $R$.  

In our generalization, if $s_i$ has an attracting cycle $P$, then $f$ will have a (nonminimal) saddle basic set 
$\Lambda(P)$ over $L_i$ with $f|\Lambda(P)$ of topological entropy $\log k_i$
(compare \cite[Proposition 3.8]{MatDil}). 
From the construction, we will see that the saddle set 
$\Lambda$ for $f$ is precisely the union of the saddle basic sets $\Lambda(P)$ 
over all attracting cycles of $s_1$ and $s_2$.


Figure~\ref{fig:cantcircbas} shows slices of $K$ for a map with $s_1(w) = w^2, s_2(w) = w^2-1$.

\begin{figure}
\begin{center}
  \drawfigcantcircbas
\end{center}
\caption{\label{fig:cantcircbas}
Let $ f(z,w)= (z^2 - 20, w^2 + z^2 - 0.9z - 20.5).$  This map appears to be in the same class as Theorem~\ref{thm:s1s2}, for 
$s_1(w) = w^2, s_2(w) = w^2-1$.
 As in Figure~\ref{fig:twbas}, we show slices of $K$.  In the center is $K_p$ (a real Cantor set).  
On top (left to right) are fibers: $z=-5, -4.99, -4.014, -4$, and on bottom (left to right) are: $z=3.998, 4, 4.886, 5$.  Note $K_{5}$ is the unit disk, and $K_{-5}$ maps onto it.  Also $K_{-4}$ is a basillica, and $K_{4}$ maps onto it. 
}
\end{figure}

\medskip\noindent{\bf Proof of Theorem \ref{thm:s1s2}.}
For clarity of exposition, we begin with a detailed 
construction for the case where $s_1(w) = s_2(w)$. Hence, fix a hyperbolic, monic polynomial $s$ of degree $d \geq 2$, and positive integers  $k_1+k_2 =d$.  

First we define a base polynomial $p = p_{k_1,r,R}$,
 for any $0 < r < R$ (the constants $r$ and $R$ will be chosen later based on $s$).  Fix distinct 
points $\xi_1, \ldots, \xi_{k_1}$ in $D_1 = D(R,r)$ and $\xi_{k_1+1}, \ldots, \xi_d$ in 
$D_2 = D(-R,r)$.  Let $p_0(z) = \prod_{j=1}^d (z-\xi_j)$.  Set 
$p(z) = ap_0(z)$ with $a>>0$ chosen so that $p^{-1}(D(0,2R))$ is a disjoint
union of $d$ disks, compactly contained in $D_1\sqcup D_2$, each univalently
mapped by $p$ onto $D(0,2R)$.   For $i = 1, 2$, define
 $$L_i = \{z\in\CC: p^n(z)\in D_i \mbox{ for all } n\geq 0\}.$$

Then $L_i$ is a forward-invariant subset of $J_p$, and $L_i$ is a Cantor set if $k_i>1$
or a single point if $k_i=1$.  In fact,
$p|L_i$ is isomorphic to the full one-sided shift on $k_i$ symbols.  Note in the case $d=2$ that
$L_1$ and $L_2$ are the two fixed points of $p$.  Finally, since $L_i \subset D_i,$ we see $L_1$ and $L_2$ are disjoint. Hence for any choice of $0 < r < R$, the map $p_{k_1, r, R}$ satisfies (1) of the theorem.
Note also that $p$ is hyperbolic with no attracting cycles.

  Define a norm $\|\cdot \|$ on the space
of polynomials of degree $d$ as the maximum of the absolute values of the coefficients. 
Choose $0< r_0 < 1$ small so that if $\{s_n\}$ is 
any sequence of polynomials with $\|s-s_n\| < 2r_0$ for all $n$, then the composition sequence
$\{s_n\circ\cdots\circ s_1\}_{n\geq 1}$ is hyperbolic with a uniform postcritical distance to the sequence
Julia sets (see \cite[Corollary 3.2]{MCrand}).  In particular, the Julia set for the composition sequence
will be a small perturbation of $J_s$.  
  
\begin{lem}  \label{lem:Mr}
There exist $M>r_0$ and $0< r< r_0$ so that 
\begin{itemize}
\item[(i)]		$M^d/18 > 2M$, 
\item	[(ii)]		$|w|>M \implies  |t(w)| \geq |w|^d/2$ for all $t$ with $\|s-t\|<2r$, 
\item[(iii)]		$\sup_{|w|\leq M} |s(w)| \leq 3M^d/2$, and
\item	[(iv)]		for any sequence $\{s_n\}$ with $\|s_n-s\| < 2r$ for all $n$, the critical 	
			points of $s$ and their 
			images under the composition sequence
		  	$s_n\circ\cdots\circ s_1$ are uniformly bounded away from 
		  	the annulus $\mathcal{A}_M = \{M\leq |w| \leq 3M\}$ and the union of 
				its preimages, $S = \bigcup_{k\geq 0} s^{-k}(\mathcal{A}_M)$.
\end{itemize}
\end{lem}

\begin{proof}
Properties (i), (ii), and (iii) can clearly be satisfied by choosing $M$ large enough and $r$
small.  Properties (i) and (ii) imply that the filled Julia set $K_s$ of $s$ is contained in the 
disk $D(0,M)$.  Therefore, 
property (iv) is only relevant if $s$ has escaping critical points, because $s$ is hyperbolic
and $S$ accumulates on the Julia set $J_s$.  Let 
$$G_s(z) = \lim_{n\to\infty} \frac{1}{d^n} \log^+|s^n(z)|$$
be the escape-rate function for $s$.   By selecting $M$ large enough, 
the images of the critical points of $s$ under iterates of $s$ 
can be arranged to be disjoint from $\mathcal{A}_M$ 
because the modulus $\operatorname{mod} \mathcal{A}_M = \log 3$ is independent of $M$, 
while a fundamental annulus $\{z\in\CC: c< G_s(z) < dc\}$ has modulus $\to \infty$ as $c\to \infty$. 
It follows that the postcritical set of $s$ is uniformly bounded away from $S$.  
Property (iv) holds for nearby sequences by continuity.
\end{proof}

Set
  $$R = 2 M^d - r,$$
and let $p = p_{k_1,r,R}$.
For our skew product map, we define 
  $$f(z,w) = (p(z), s(w) + p(z)-z).$$

\begin{lem} \label{lem:Kbounds}
The fiber filled Julia sets for $f$ satisfy
$$K_z \subseteq \left\{ \begin{array}{ll} 
		D(0,3M) & \mbox{ for all } z\in J_p  \\
		\mathcal{A}_M = \{M\leq |w| \leq 3M\} & \mbox{ for all }  z\in D_i\cap p^{-1}D_j\cap J_p, i\not= j \\
		D(0,M) & \mbox{ for all } z\in L_1\cup L_2
		\end{array} \right.  $$  
\end{lem}

\begin{proof}
From the definition of $f$ and the choice of $M$ and $r$ in Lemma \ref{lem:Mr}, we have
that $|w| \geq 3M$ and $z\in J_p$ implies: 
\begin{eqnarray*}
|q_z(w)| = |s(w) + p(z)-z|
		&\geq& ||s(w)| - |p(z)-z|| \geq \frac12 |w|^d - 2(R+r) \\
		&=& \frac12 |w|^d - 4 M^d = |w|^d(\frac12 - 4(M/|w|)^d)\\
		&\geq&  |w|^d/18 > 2|w|,
\end{eqnarray*}
and for $z\in D_i\cap p^{-1}D_j, i\not=j$,  the inequality $|w| \leq M$ implies:
\begin{eqnarray*}
|q_z(w)| &\geq& ||s(w)| - |p(z)-z|| \geq 2(R-r) - \frac{3}{2}M^d  \\
		&=&  4 M^d -4r - \frac{3}{2}M^d > 45 M - 4r > 41M.
\end{eqnarray*}
Therefore these points escape to infinity under iteration of $f$.  

For $z\in L_1\cup L_2$, the sequence of polynomials $\{q_{p^n(z)}\}$ satisfies
$\|s -q_{p^n(z)}\|< 2r$ for all $n\geq 0$;  by the choice of $r$,  properties (i) and (ii)
of Lemma \ref{lem:Mr} imply that $K_z\subset D(0,M)$.  
\end{proof}

For $z\in L_1\cup L_2$, the composition sequences 
$\{Q_z^n\}$
are hyperbolic by our choice of $r$ with uniform postcritical distance
\cite[Theorem 1.3]{MCrand}, so $f$ is vertically 
expanding over $L_1\cup L_2$.  

The critical points $C_z$ of $f$ over $z\in J_p$ coincide with the critical points of $s$.  
For $z\in J_p \setminus (L_1 \cup L_2)$, let $z_n = p^n(z)$.  There is a smallest integer $N$ so that $z_0, \dots, z_N \in D_i$ and $z_{N+1} \in D_j$ (with $j\not=i$).  
From (iv) of Lemma~\ref{lem:Mr},  the images of the critical points $Q_z^n(C_z)$ remain in $\{|w| < M\}\cup\{|w| >3M\}$ for all $n\leq N$.  Furthermore, $K_{z_N} \subseteq \mathcal{A}_M$ by Lemma \ref{lem:Kbounds}. 
Let $S = \bigcup_{k\geq 0} s^{-k}(\mathcal{A}_M)$.  
From Lemma \ref{lem:Mr}, the postcritical set of any sequence $\{s_n\}$ with $\|s - s_n\|<2r$ for all $n$ is uniformly bounded away from $S$.  By invariance of the filled Julia sets, the postcritical points $Q_z^n(C_z)$ are therefore uniformly bounded away from the filled Julia set $K_{z_n}$ for all $n\leq N$.  
For $n\geq N+1$, we have $|Q_z^n(C_z)| > 6M$, so these points are also uniformly
bounded away from $K_{z_n}$.  

We conclude that $f$ is vertically expanding over all of $J_p$. 
Since $p$ has no attracting periodic points, $f$ is Axiom A by Theorem~\ref{thm:Jthm8.2}.
 Note that all critical points
over $J_p\setminus (L_1\cup L_2)$ escape. Combining this with Proposition~\ref{prop:Jprop2.3} yields (2) of the theorem.

\bigskip
In the case of distinct $s_1$ and $s_2$, each monic, hyperbolic 
polynomials of degree $d$, choose $r$ and $M$ as in Lemma~\ref{lem:Mr} to work for both $s_1$ and 
$s_2$, set $R = 2M^d-r$ and $p = p_{k_1, r/2, R}$.
Write
  $$s_1(w) = w^d + t_1(w)$$
  $$s_2(w) = w^d + t_2(w)$$
where $t_i(w) = O(w^{d-1})$ for $i = 1,2$.  Set 
  $$L(z,w) = w^d + \frac{z+R}{2R} t_1(w) - \frac{z-R}{2R} t_2(w),$$
and
  $$f(z,w) = (p(z), L(z,w) + p(z)-z).$$
Then the arguments above show that $f$ is vertically expanding over $J_p$ (hence Axiom A), and that (1) and (2) hold, and further, $f$ behaves as a 
small perturbation of the product $(p|L_1) \times s_1$ over $L_1$ and as a small
perturbation of $(p|L_2)\times s_2$ over $L_2$, which establishes (3) of the theorem.

\bigskip

For (4) of the theorem, first recall  that $\Lambda = \Apt$ (Lemma~\ref{lem:AptisL}), and note $\Apt \subset (L_1 \cup L_2) \times \CC$, since we showed above that all critical points over $J_p \setminus (L_1 \cup L_2)$ escape.  Since $L_1$ and $L_2$ are disjoint and each is forward invariant,
 we conclude $\Lambda = \Lambda_1 \cup \Lambda_2$ with $\Lambda_i \subset (L_i \times \CC)$, and any basic set in $\Lambda$ is contained in one of $\Lambda_1$ or $\Lambda_2$.  Note also $\Lambda_i$ is the closure of the saddle periodic points of $f$ in $L_i \times \CC$.

Let $\alpha_i$ be any periodic point of $p$ in $L_i$, say of period $n$ (at least one $\alpha_i$ exists since $p|L_i$ is the full one-sided shift on $k_i$ symbols). Then by (3), 
$f^n|_{ \{ \alpha_i \} \times \CC} = Q^n_{\alpha_i} \colon \{ \alpha_i \} \times \CC \to \{ \alpha_i \} \times \CC$ is a small perturbation of $s_i^n$, so its attracting cycles are perturbations of those for $s_i^n$.  Thus the saddle basic sets of $f$ in $L_i \times \CC$ are in one-to-one correspondence with attracting periodic points of $s_i$.
Hence $\Lambda_i \neq \emptyset$ precisely when $s_i$ has an attracting cycle, which, 
since $s_i$ is hyperbolic, is equivalent to $s_i$ having a critical point which does not escape.



\bigskip
To show (5) of the theorem,
assume that $s_i$ has an attracting cycle $P = \{ w_1, \ldots, w_m \} \subset\CC$.
Then as above, $f$ has an associated saddle basic set $\Lambda(P)$ over $L_i$.
Let $c$ be a critical point of $s_i$ such that $s_i^n(c)\to P$ as $n\to\infty$.

Fix $z\in J_p\cap D_i$ and let $z_{-k}\in p^{-k}(z)$ be a sequence of 
preimages such that $z_{-k}\in D_i$ for all $k$.  By the construction of $f$
(in particular, the choice of $r$), the point $f^k(z_{-k},c)$ lies in a small 
neighborhood $U_P$ of $P$ in $\{z\}\times\CC$ for all $k$ sufficiently large.  
Therefore this neighborhood contains a point
in the accumulation set $\A$.  Consequently, for all $z\in J_p$, we have
  $$\A\cap (\{z\}\times U_P) \not= \emptyset.$$   

On the other hand, $\Acc = \Apt$ because the connected components of $J_p$
are points, while $\Apt = \Lambda$ by Lemma~\ref{lem:AptisL}.  As $\Lambda
\subset (L_1\cup L_2)\times\CC$, we conclude that $\A \not= \Acc$.   

Finally, recall that products must have equality of $\Lambda = W^u(\Lambda) \cap (J_p \times \CC)$
by Proposition \ref{prop:products}, and this equality is preserved in hyperbolic
components (Proposition~\ref{prop:Lcomponent}).  Therefore, in the case
that $\Lambda\not=\emptyset$, we can conclude that $f$ is not in the same 
hyperbolic component as a product.  This concludes the proof of Theorem \ref{thm:s1s2}.
\qed

\begin{rem}
The extension of an Axiom A $f$ (as constructed in Theorem \ref{thm:s1s2}) 
to the line at infinity is
$\zeta \mapsto \zeta^d + a$, where $a$ was chosen in defining $p$ at the 
beginning of the proof, and $a$ is very large (so that
$p(z) = a p_0(z)$ has a Cantor Julia set).  Thus $f$ is expanding on
the Cantor Julia set on the line at infinity, so $f$ is Axiom A on $\Pt$.
Thus $f$ satisfies (4) of Theorem~\ref{thm:examples}.
\end{rem}


\section{Remaining Questions}

\begin{question}  \label{question1}
In our examples, we focused on maps of degree two.  In higher degree, more varied phenomena than we discussed might occur.  

Let $f(z,w) = (p(z), q(z,w))$ be an Axiom A polynomial skew product of degree $d$.
In the space of polynomial maps of $\CC$ of degree $d$,
let $\mathcal{E}_d$ be the hyperbolic polynomials with Cantor Julia set
and $\mathcal{HC}_d$ be the hyperbolic polynomials with connected Julia set.
If $p \in \mathcal{E}_d$, then $\Apt = \Acc$, and if $p \in \mathcal{HC}_d$, then $\Acc = \A$.  All hyperbolic polynomials of degree two are either in $\mathcal{E}_2$ or $\mathcal{HC}_2$.  But in higher 
degree this is not the case.   

Thus we ask: do there exist Axiom A polynomial skew products $f_1, f_2, f_3$ such that $p_k \not\in (\mathcal{E}_d \cup \mathcal{HC}_d), k=1,2,3$ and:
\begin{enumerate}
\item	$\Apt = \Acc \not= \A$ for $f_1$;
\item	$\Apt \not= \Acc = \A$ for $f_2$;
\item	$\Apt\not= \Acc \not= \A$ for $f_3$?
\end{enumerate}
\end{question}

\begin{question}  \label{question2}
Does there exist a characterization of the equality $\Acc = \A$ 
in a similar spirit to Lemma~\ref{lem:AptAccchar} or Theorem \ref{thm:Lbigequiv}?  And is this equality 
preserved in hyperbolic components?  

Propositions~\ref{prop:Lcomponent}
and~\ref{prop:AptisAcc} imply that (1) of Question \ref{question1} is preserved in hyperbolic components.   We are asking if the same is true for (2) and (3).
\end{question}

\begin{question}  \label{question3}
Nekrashevych \cite{Nek} shows that the \textit{rational} skew product of $\Pt$ given by
$$
R(z,w) = \(1 - 1/z^2, 1 - w^2/z^2 \)
$$
is Axiom A, with connected base Julia set, and all fiber Julia sets connected, but such that not all fibers are homeomorphic (for example, over the fixed points of the base map, one fiber map is the rabbit, while another one is the airplane).   This suggests there is no dynamical obstruction to a polynomial skew product of $\Ct$ which is fully connected yet with varying fiber dynamics, but no such example has been exhibited.

By Corollary~\ref{cor:productcomponent}, such a map would not be in the same hyperbolic component as any product, and by Corollary~\ref{cor:fconnected}, it would satisfy $\Apt = \Acc = \A$.
\end{question}

\bibliographystyle{alpha}
\bibliography{Skew}

 \end{document}